\documentclass[10pt]{article}
\usepackage{amssymb}
\usepackage{amsmath}
\usepackage{multirow}
\usepackage{graphicx}

\newtheorem{theorem}{Theorem}[section]

\newtheorem{proposition}[theorem]{Proposition}

\def\e{\mathrm e}

\usepackage{xcolor}

\xdefinecolor{gray}{rgb}{0.4,0.4,0.4}
\xdefinecolor{blue}{RGB}{58,95,205}
\definecolor{dkgreen}{rgb}{0,0.6,0}
\definecolor{gray}{rgb}{0.5,0.5,0.5}
\definecolor{mauve}{rgb}{0.58,0,0.82}
\definecolor{darkblue}{rgb}{0.0,0.0,0.6}
\definecolor{cyan}{rgb}{0.0,0.6,0.6}
\usepackage{listings}

\graphicspath{{.}{Figures/}}

\begin{document}

\title{A tutorial on estimator averaging in  spatial point process models}

\author{F. Lavancier and P. Rochet}

\maketitle

\begin{abstract}
Assume that several competing methods are available to estimate a parameter in a given statistical model.
The aim of estimator averaging is to provide a new estimator,  built as a linear combination of the initial estimators, that achieves better properties, under the quadratic loss,  than each individual initial estimator. This contribution provides an accessible and clear overview of the method, and investigates its performances on standard spatial point process models. It is demonstrated that the average estimator clearly improves on standard procedures for the considered models. For each example, the code to implement the method with the \texttt{R} software (which only consists of  few lines) is provided.
\end{abstract}

\begin{center}
\end{center}

\section{Introduction}\label{sec:intro}

Assume one needs to estimate a real parameter $\theta$ when several possibly competing methods are known to the statistician, leading to a collection of estimators $\hat \theta_1,...,\hat \theta_J$ with $J \geq 2$. As it has been observed in numerous practical situations, the initial estimators $\hat \theta_1,...,\hat \theta_J$ may contain complementary information on the parameter of interest so that choosing a single one in the collection might not be optimal. A well-spread idea in this situation is to consider a linear combination of the $\hat \theta_j$'s that would hopefully preserve each individual quality. A final estimator is then sought as a combination
$$ \hat \theta_\lambda = \sum_{j=1}^J \lambda_j \hat \theta_j \ \ \text{ subject to } \ \ \sum_{j=1}^J \lambda_j = 1,  $$ 
where $\lambda = (\lambda_1,...,\lambda_J)^\top \in \mathbb R^J$ is the vector of weights that has to be estimated. The main purpose is to provide an estimator that would  perform at least as well as the best estimator in the initial collection, or even better if possible.

The gains of considering combinations of estimators are well established in the literature in the particular case of predictors in regression models as in \cite{hansen2007least,hjort2003frequentist,MR2351101,yang2004aggregating} or forecasts in time series, see  \cite{bates1969combination,elliott2011averaging,timmermann2006forecast}. Depending on the community, these techniques are commonly referred to aggregation or model averaging. Recently, a general methodology for estimator averaging was proposed in \cite{lavancier2016general} in an attempt to extend model averaging beyond prediction purposes. In the present paper, we give a clear review of this method, providing the code in the \texttt{R} software for its implementation, and we investigate the efficiency of averaging in some common spatial  point process models. Specifically we consider  the estimation of the intensity  of an inhomogeneous spatial Poisson point process, the estimation of the parameters in a determinantal point process (a model for regular point patterns), in a Thomas process (a model for clustered point patterns), and in a Boolean model (the basic model for random sets). We argue that the averaging procedure is particularly well suited to these models due to the lack, in most cases, of a universal single best estimation method.  For these examples, we demonstrate that  the average estimator performs better (in the mean square sense) than the best initial estimator, which conveys that none of the current methods is able to gather the whole information available from the data. Moreover, the averaging procedure allows not only to get a better estimate, but also provides for the same price an estimation of its mean square error, which allows to construct confidence intervals without further effort.

For each example,  we describe the full implementation of  the  averaging procedure in the software \texttt{R} \cite{R}. We point out that the code takes only few lines of scripts that mainly rely on routines available in the package \texttt{spatstat} \cite{baddeley:rubak:turner:15, A-BadTur05}. We also indicate the  CPU time on a laptop equipped with an Intel i7, 2.6GHz processor.  It varies from 20 seconds to 3 minutes without parallelization, that can become few seconds if a parallelization procedure is used in the resampling step of the method. Depending on the considered model, the CPU time  is either due to the simulation of the model or to the estimating procedures, that both have to be repeated several times in the resampling step. 

The rest of the paper is  organized as follows. Section~\ref{sec:method} contains an accessible overview of the averaging procedure, emphasizing the key choices for its implementation and providing the associated code in \texttt{R}. In Section~\ref{sec:appli}, we 
apply this method to the models of spatial statistics listed above, demonstrating the relevance of the procedure in these cases. 
We conclude in Section~\ref{sec:conclusion} with a brief summary of our study and some general recommendations.

\section{Description of the method}\label{sec:method}

This section is devoted to the description of the averaging procedure introduced in \cite{lavancier2016general}. For ease of understanding, the framework is simplified so as to fit more precisely with the spatial statistical models studied in this paper. Nevertheless, the method described here remains general and suitable for most parametric and semi-parametric models, whether they concern spatial statistics or more generic frameworks.

\subsection{Averaging of a real-valued parameter}\label{sec:real}

Let $\theta \in \mathbb R$ and $\hat \theta_1,...,\hat \theta_J$ a collection of $J \geq 2$ estimators of $\theta$. For $\lambda=(\lambda_1,...,\lambda_J)^\top \in \mathbb R^J$ a vector of weights such that $\sum_{j=1}^J \lambda_j = 1$, we are interested in the performance of the average estimator
$$ \hat \theta_\lambda = \sum_{j=1}^J \lambda_j \hat \theta_j .$$
The condition $\sum_{j=1}^J \lambda_j = 1$ was originally proposed in \cite{bates1969combination} as a way to preserve the unbiasedness (if so) of the initial estimators, but may in fact be important for deeper reasons. Actually, it is argued in \cite{lavancier2016general} that the condition $\sum_{j=1}^J \lambda_j = 1$ is crucial when the parameter $\theta$ is real-valued, as in the present situation. However, this restriction can be safely overlooked for more complex parameters, as in the 
 context of non-parametric regression, see \cite{MR2351101}. 

If all the initial estimators $\hat \theta_j$ are square integrable, the combination $\lambda^*$ minimizing the quadratic risk expresses easily in function of the mean-square error (MSE) matrix $\Sigma$ with general term $\Sigma_{ij} = \mathbb E \big[(\hat \theta_i - \theta)(\hat \theta_j- \theta)\big]$, $i,j=1,...,J$. Indeed, due to the condition $\sum_{j=1}^J \lambda_j = 1$, the quadratic error of $\hat \theta_\lambda$ simplifies into
\begin{equation}\label{mse}
 \mathbb E \big(\hat \theta_\lambda - \theta \big)^2 = \mathbb E \Big( \sum_{j=1}^J \lambda_j (\hat \theta_j - \theta) \Big)^2 = \lambda^\top \Sigma \lambda\end{equation}
where $\lambda^\top$ denotes the transpose of $\lambda$. 
Thus, the expression of the best linear combination $\lambda^*$ which minimizes the MSE and determines the so-called \textit{oracle} $\hat \theta^* := \sum_{j=1}^J \lambda^{*}_j \hat \theta_j$ follows as the solution of a simple constrained optimization problem
\begin{equation}\label{lambda_one}
 \lambda^* = \arg \min_{ \lambda \in \mathbb R^J: \lambda^\top \mathbf 1 =1} \lambda^\top \Sigma \lambda = \frac{\Sigma^{-1} \mathbf 1} {\mathbf 1^\top \Sigma^{-1} \mathbf 1}, 
 \end{equation}
where $\mathbf 1 = (1,...,1)^\top \in \mathbb R^J$. The oracle is unknown in practice, but can be approximated whenever an estimate $\hat \Sigma$ of the MSE matrix is available. The average estimator $\hat \theta$ is then constructed as an approximation of the oracle obtained by replacing $\Sigma$ by $\hat \Sigma$ in the expression of $\lambda^*$, that is
\begin{equation}\label{weight} 
 \hat \theta = \hat \theta_{\hat \lambda}  \quad\text{with}\quad \hat \lambda = \frac{\hat \Sigma^{-1} \mathbf 1} {\mathbf 1^\top \hat \Sigma^{-1} \mathbf 1},  \end{equation}
provided that $\hat \Sigma$ is non-singular. Remark that the computational cost of the method comes essentially from producing the matrix $\hat \Sigma$. Once $\hat \Sigma$ (denoted below by  \texttt{hatSigma}) is available, deducing the weights $\hat \lambda$ in \texttt{R} is straightforward:
\begin{lstlisting}
invhatSigma <-solve(hatSigma)
weights <-rowSums(invhatSigma)/sum(invhatSigma)
\end{lstlisting}
An additional benefit of this method is to provide an estimation of the mean square error of the resulting estimator $\hat\theta$ without further effort. Indeed, under some conditions the MSE of $\hat\theta$ becomes asymptotically equivalent to  the MSE of the oracle (see below and \cite{lavancier2016general}), which is given by ${\lambda^*}^\top \Sigma \lambda^*=(\mathbf 1^\top \Sigma^{-1} \mathbf 1)^{-1}$ and can naturally be estimated by 
\begin{equation}\label{MSEAV}\widehat{MSE}(\hat\theta)=(\mathbf 1^\top \hat \Sigma^{-1} \mathbf 1)^{-1}.\end{equation}
In  \texttt{R}, this is simply:
\begin{lstlisting}
MSE_AV <-1/sum(invhatSigma)
\end{lstlisting}

 The estimation of $\Sigma$ can be carried out with the same data as those used to produce the initial estimators $\hat \theta_1,...,\hat \theta_J$. In particular, the averaging procedure does not require the independence between $\hat \Sigma$ and the $\hat \theta_j$'s. For parametric models, as this is commonly the case in spatial statistics, $\Sigma$ can simply be estimated by   a parametric bootstrap procedure, specifically: 
\begin{enumerate}
\item Choose an initial consistent estimate $\hat\theta_0$ (typically one of the initial estimators, or their simple average, see also the discussion in Section~\ref{sec:complement}).
\item Simulate $N$ samples according to the model with the previous estimate as a parameter.
\item For each sample $b=1,\dots,N$, compute the estimates $\hat \theta_1^{(b)},\dots,\hat \theta_J^{(b)}$ where the superscript $(b)$ emphasizes the dependence on the sample $b$.
\item Deduce an estimation of each entry of $\Sigma$ as $\hat \Sigma_{ij}=1/N \sum_{b=1}^N (\hat\theta_i^{(b)} - \hat\theta_0)(\hat\theta_j^{(b)} - \hat\theta_0)$.
\end{enumerate}
This procedure is used with $N=100$ in all our examples in  Section~\ref{sec:appli}, where we provide the associated  \texttt{R} code. Of course the larger $N$, the better the approximation in the fourth step. But it turns out that $N=100$ appeared as a good compromise in our examples to get a fast and decent  approximation. Alternative methods to estimate $\hat\Sigma$, in particular for semi or non-parametric models, are discussed in Section~\ref{sec:complement} and in \cite{lavancier2016general}.\\

From a theoretical point of view, the performance of the average estimator $\hat \theta$ can be measured by how well $\hat \Sigma \Sigma^{-1}$ approximates the identity matrix. A non-asymptotic bound on the distance to the oracle is derived in Theorem 3.1 in \cite{lavancier2016general}, although the behavior of the error term is hardly tractable in practice. This result guarantees nevertheless the asymptotic optimality of the average estimator in the following sense. 

\begin{proposition}\textbf{\emph{\cite{lavancier2016general}}} \label{prop} If $\ \hat \Sigma \Sigma^{-1}$ converges in probability to the identity matrix when the sample size tends to infinity, then
$$ \hat \theta - \theta = \hat \theta^* - \theta + o_p \big( \mathbb E (\hat \theta^* - \theta)^2 \big).$$
\end{proposition}

Note that the crucial condition on $\ \hat \Sigma \Sigma^{-1}$ above holds true if $\hat\Sigma$ is obtained by parametric bootstrap, provided $\Sigma$ is a sufficiently smooth function of the parameters.
This is the case for all parametric models considered in  Section~\ref{sec:appli}.
Under additional technical assumptions on the moments of $\hat \Sigma$ and the $\hat \theta_j$'s (see \cite{lavancier2016general}), one can deduce the asymptotic optimality in $\mathbb L^2$

$$  \mathbb E (\hat \theta - \theta)^2 = \big( 1+ o(1) \big) \ \mathbb E(\hat \theta^* - \theta)^2.  $$

\subsection{Averaging with foreign estimators}\label{sec:foreign}

Another important advantage of the averaging procedure is that it allows to use information contained in estimators of other parameters. Assume that the true distribution of the observation depends on both $\theta$ and a nuisance parameter $\eta \in \mathbb R$ with a collection of estimators $\hat \eta_1,..., \hat \eta_K$ also available for $\eta$. We investigate situations where the use of the $\hat \eta_k$'s can improve the estimation of $\theta$. In this context, the $\hat \eta_k$'s are referred to as foreign estimators. 

The decision to include foreign estimators is in particular motivated by the relative efficiency of the estimations. For instance, if the parameter $\theta$ is known to be poorly estimated, using another better estimated parameter $\eta$ generally tends to improve the performance of  the $\hat \theta_j$'s, if the $\hat \theta_j$'s and the  $\hat \eta_k$'s are correlated. 

\bigskip

\textit{Remark.} We consider only a one real-valued nuisance parameter $\eta$  for simplicity. This is the framework for the estimation  of the Boolean model treated in Section~\ref{sec:bool}. Nevertheless, the method can be easily extended to situations with several nuisance parameters (see \cite{lavancier2016general} for more details). In Section~\ref{sec:Thomas}, we apply it to the estimation of three parameters, allowing for the inclusion of two foreign parameters for each. \bigskip

The foreign estimators $\hat \eta_k$'s are included in the estimation of $\theta$ by considering an additional linear combination of the $\hat \eta_k$'s with the weights summing to zero. Thus, a final estimate of $\theta$ is sought as a combination of the $\hat \theta_j$'s and $\hat \eta_k$'s,
$$ \hat \theta_{\lambda,\mu} = \sum_{j=1}^J \lambda_j \hat \theta_j + \sum_{k=1}^K \mu_k \hat \eta_k \ \  \text{ subject to } \ \ \sum_{j=1}^J \lambda_j = 1 \ \ \text{and} \ \  \sum_{k=1}^K \mu_k = 0.    $$
The main reason for imposing the $\mu_k$'s to sum to zero is that the oracle $(\lambda^*,\mu^*)$ can still be expressed in function of the MSE matrix, but this time of the whole collection $(\hat \theta_1,...,\hat \theta_J,\hat \eta_1,...,\hat \eta_K)$. To see this, write
\begin{equation}\label{eqmforeign} 
\mathbb E \big(\hat \theta_{\lambda,\mu} - \theta \big)^2  = \mathbb E \Big(\sum_{j=1}^J \lambda_j (\hat \theta_j - \theta) + \sum_{k=1}^K \mu_k (\hat \eta_k - \eta) \Big)^2 = (\lambda^\top,\mu^\top) \ \Sigma \ \Big(\begin{array}{c} \!\!\lambda\!\! \\ \!\!\mu\!\! \end{array} \Big)
\end{equation}
where $\Sigma$ designates here the $(J+K)\times (J+K)$ MSE matrix of $(\hat \theta_1,...,\hat \theta_J,\hat \eta_1,...,\hat \eta_K)$. Here again, the optimal combination minimizing \eqref{eqmforeign} subject to $\sum_{j=1}^J \lambda_j = 1$ and $\sum_{k=1}^K \mu_k = 0$ has a simple expression in function of $\Sigma$, described below. 

For sake of completeness, let us consider the full problem of estimating both $\theta$ and $\eta$, using foreign estimators for each. 
This means that we seek the average estimators  $\hat \theta_{\lambda,\mu}$ and  $\hat \mu_{\lambda',\mu'}$, built as explained above with the constraints $\sum_{j=1}^J \lambda_j = 1$, $\sum_{k=1}^K \mu_k = 0$ for $\hat \theta_{\lambda,\mu}$, and $\sum_{j=1}^J \lambda'_j = 0$, $\sum_{k=1}^K \mu'_k = 1$ for $\hat \mu_{\lambda',\mu'}$.  As proved in \cite{lavancier2016general}, the optimal weights $(\lambda^*,\mu^*)$ and $({\lambda'}^*,{\mu'}^*)$ minimizing respectively $\mathbb E \big(\hat \theta_{\lambda,\mu} - \theta \big)^2$ and  $\mathbb E \big(\hat \eta_{\lambda',\mu'} - \eta \big)^2$
are given by the $(J+K)\times 2$ matrix
\begin{equation}\label{oracle_max} 
 \left( \begin{array}{cc} \lambda^* & {\lambda'}^* \\ \mu^* & {\mu'}^*\end{array} \right) = \Sigma^{-1} \operatorname{L} (\operatorname{L}^\top  \Sigma^{-1} \operatorname{L})^{-1}
\end{equation}
where $\operatorname{L}$ denotes the $(J+K)\times 2$ matrix
$$\operatorname{L} = \left( \begin{array}{cc} \!\! \mathbf 1_{J} \!\! & 0  \\ 0 &  \!\! \mathbf 1_{K} \!\!   \end{array} \right), $$
with  $\mathbf 1_J = (1,...,1)^\top \in \mathbb R^J$ and $\mathbf 1_K = (1,...,1)^\top \in \mathbb R^K$. 
This solution is approximated in practice  using an estimate $\hat \Sigma$, typically obtained by parametric bootstrap as detailed in the previous section. The code in \texttt{R} to get an estimate of the matrix of optimal weights \eqref{oracle_max}, given $\hat\Sigma$ (\texttt{hatSigma}), $J$ and $K$ is as follows
\begin{lstlisting}
invhatSigma <-solve(hatSigma)
matL <-matrix(c(rep(1,J),rep(0,K),rep(0,J),rep(1,K)),ncol=2)
weights <-invhatSigma%*%matL%*%solve(t(matL)%*%invhatSigma%*%matL)
\end{lstlisting}
Here again, an estimation of the MSE matrix of $(\hat\theta,\hat\eta)$ is straightforward from~\eqref{eqmforeign} 
\begin{lstlisting}
MSE_AV <-t(weights)%*%hatSigma%*%weights
\end{lstlisting}

Remark that the optimal weights $(\lambda^*,\mu^*)$  for the estimation of $\theta$ are derived simultaneously in \eqref{oracle_max}, so that $\lambda^*$ is not equal to the optimal combination \eqref{lambda_one} when the foreign estimators are not included. In theory, including foreign estimators leads to a better oracle.  In practice, while introducing foreign estimators does produce convincing results in specific situations (see e.g. Section~\ref{sec:bool}), the simulation study suggests that it only occasionally improves on the situation where only the estimators of $\theta$ are involved in the averaging process (see Section~\ref{sec:Thomas}). However, it never dramatically reduces the performance of the average estimator. The key point is of course the quality of estimation of $\Sigma$. The question of including or not foreign estimators in the averaging process must be subject to a preliminary analysis, proper to each model,  in order to optimize the performances of the averaging procedure.

\subsection{Complementary aspects}\label{sec:complement}

The averaging procedure described in the previous sections provides a simple way to derive a single accurate solution for statistical inference when several competing methods are available. The only two factors to calibrate are the estimation of $\Sigma$ and the choice to include or not foreign estimators (the latter only if the model contains several parameters). 
We discuss in this section  the estimation of $\Sigma$ and the possibility  to account for additional constraints on the weights.\\

Concerning the estimation of $\Sigma$, we systematically use in the examples of this paper a parametric bootstrap procedure, as described in Section~\ref{sec:real}.  In a parametric model, $\Sigma$ indeed depends on the parameters and $\hat \Sigma$ can then be obtained by plug-in or, if an analytical form of $\Sigma$ in function of the parameters is not available (which is commonly the case as for our examples), by parametric bootstrap. In this situation, the performance of the average estimator can be highly dictated by the choice of the initial estimator $\hat\theta_0$ used to perform the bootstrap procedure. 
As a general recommendation, we suggest for $\hat\theta_0$ to use in that order: 1) the average of the initial estimators if they are comparable in efficiency, 2) the best overall estimator in the collection if it is known and 3) a robust estimator if the performances of the initial estimators  are very variable  depending on the true distribution of the data. 

In a semi or non-parametric model where the expression of $\Sigma$ is more complicated, other methods to estimate $\Sigma$ can be considered.  A first alternative is to use an asymptotic approximation (if available) which may lead to a simplified form of the MSE matrix, typically a parametric expression, thus easier to approximate. Of course, the asymptotic form of $\Sigma$ works all the more that the amount of data is large. A second alternative is to use standard (non-parametric) bootstrap, i.e. from random sampling on the original dataset. This solution generally well applies in situations where the data are independent and identically distributed, but is however rarely suited for spatial statistics models. We refer to \cite{lavancier2016general} for examples where these methods are applied.\\

Concerning the weights of averaging, in addition to the normalization $\sum_{j=1}^J \lambda_j = 1$ (and $\sum_{k=1}^K \mu_k = 0$ for foreign estimators) considered in the previous sections, it is possible to impose additional conditions. A natural option is to restrict to positive weights $\lambda_j$ aiming for a convex combination of the initial estimators $\hat \theta_j$. This is a natural way to get a more stable final estimator since the weights are then restricted to the interval $[0,1]$. Convex averaging may lead to a sparse combination, i.e. a solution that only involves a subset of the initial estimators, which allows to perform an indirect selection among the $\hat \theta_j$'s. Another desirable property of convex averaging arises when the parameter of interest $\theta$ lies in a convex subset of $\mathbb R$ (e.g. $\theta \geq 0$ or $\theta \in [0,1]$). In this case, the solution is guaranteed to remain in the same space as the initial estimators due to its stability by convex transformations.
As to the implementation of convex averaging, the problem of minimizing the MSE \eqref{mse}  subject to the constraints  $\sum_{j=1}^J \lambda_j = 1$ and $\lambda_j\geq 0$ has no explicit solution. It is however an easy quadratic optimization problem that can be numerically solved in \texttt{R} as follows using the package \texttt{quadprog} \cite{quadprog}.
\begin{lstlisting}
temp <-solve.QP(hatSigma,rep(0,J),cbind(rep(1,J),diag(J)),c(1,rep(0,J)),meq=1)
weights_convex <-temp$solution
\end{lstlisting}
It is worth emphasizing that the additional constraints $\lambda_j\geq 0$ in convex averaging result in a more accessible but less accurate oracle. The same remark holds for any additional constraint on the weights. Thus,  if $\Sigma$ can be suitably estimated, it is generally not too risky to consider the minimal constraint  $\sum_{j=1}^J \lambda_j = 1$ thus aiming for the best possible oracle. On the other hand, if the estimator $\hat \Sigma$ is not reliable, additional constraints on the weights can be set in order to make the oracle easier to approximate. Convex averaging is an option. Another option, in presence of many initial estimators (leading to a matrix $\Sigma$ difficult to estimate) is to consider the combination of a restricted number of estimators, typically two or three. This can be achieved by introducing the constraint that at least $J-2$ (or $J-3$) weights must be zero. Some preliminary simulations not shown in this paper suggest that this strategy is promising. However, due to its non-convexity, this setting is no longer covered by the theoretical guarantees provided in  \cite{lavancier2016general}. In the specific framework of Gaussian regression with sample splitting, this so-called subset aggregation strategy is investigated in  \cite{MR2351101}, see also \cite{gaiffas}. 
The study of this strategy in a more general setting, as in the present paper, is the subject of a work currently in progress.

\section{Application to spatial statistics models}\label{sec:appli}

\subsection{Inhomogeneous Poisson point process}
\begin{figure}
\centering
  \includegraphics[scale=.14]{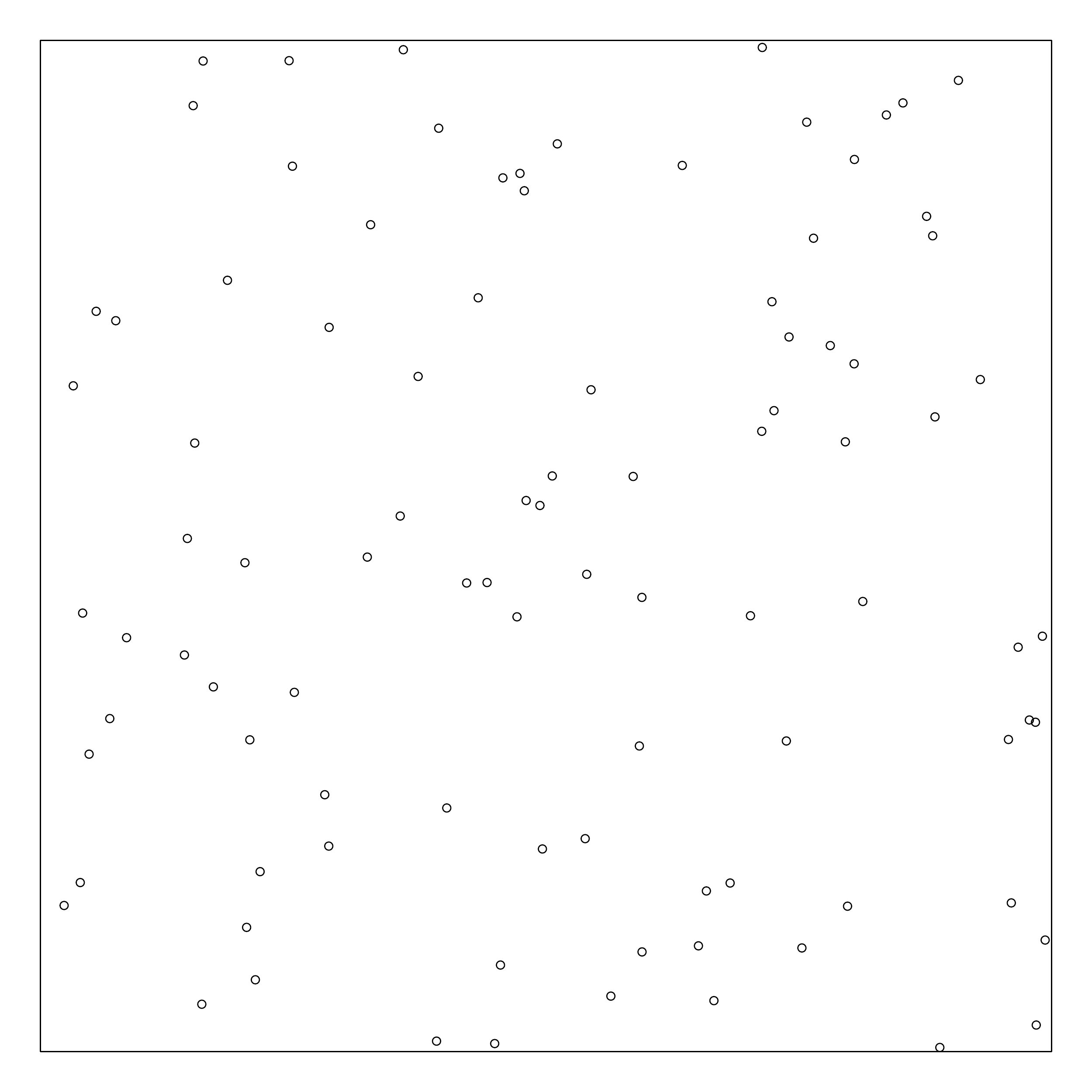}
  \includegraphics[scale=.14]{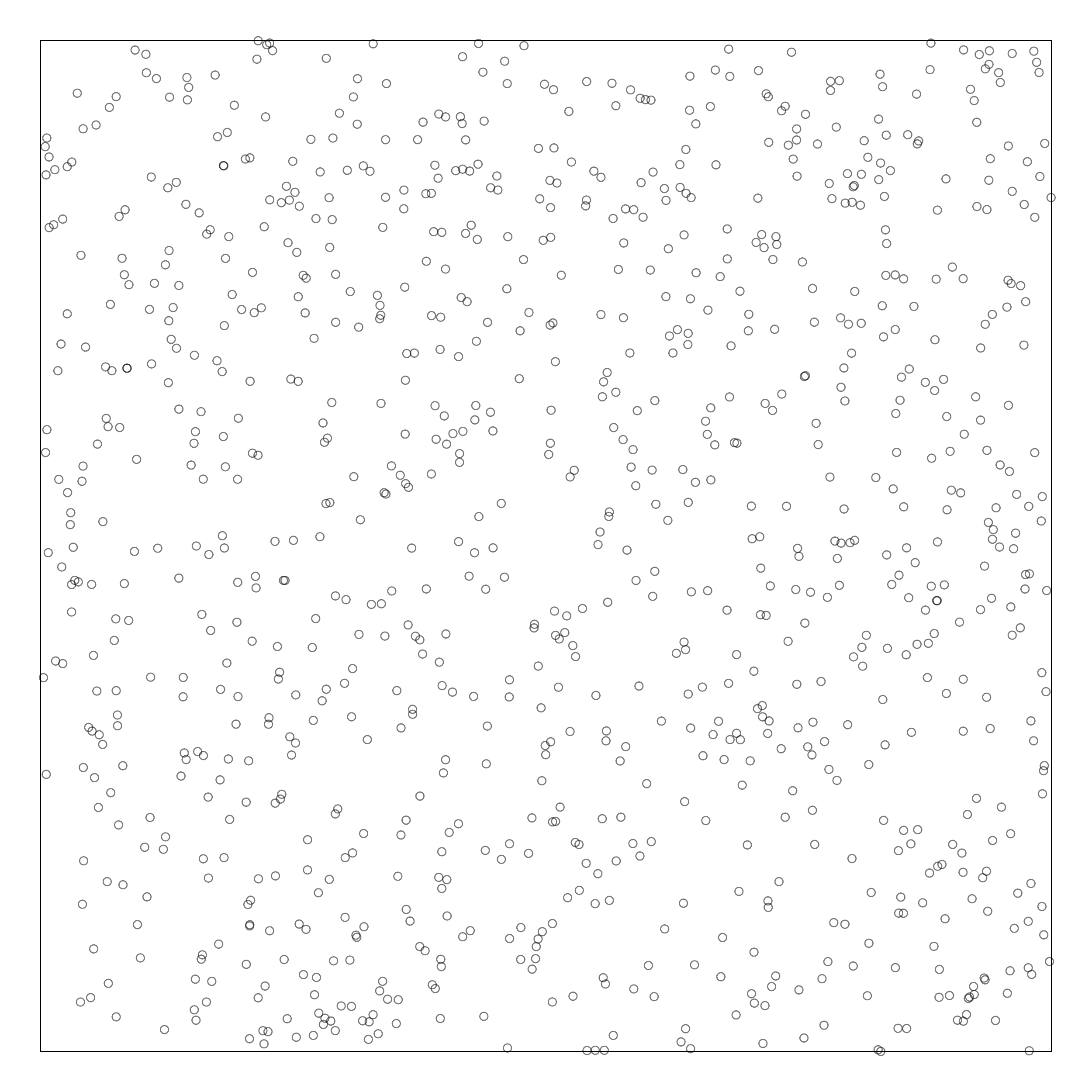}
  \includegraphics[scale=.14]{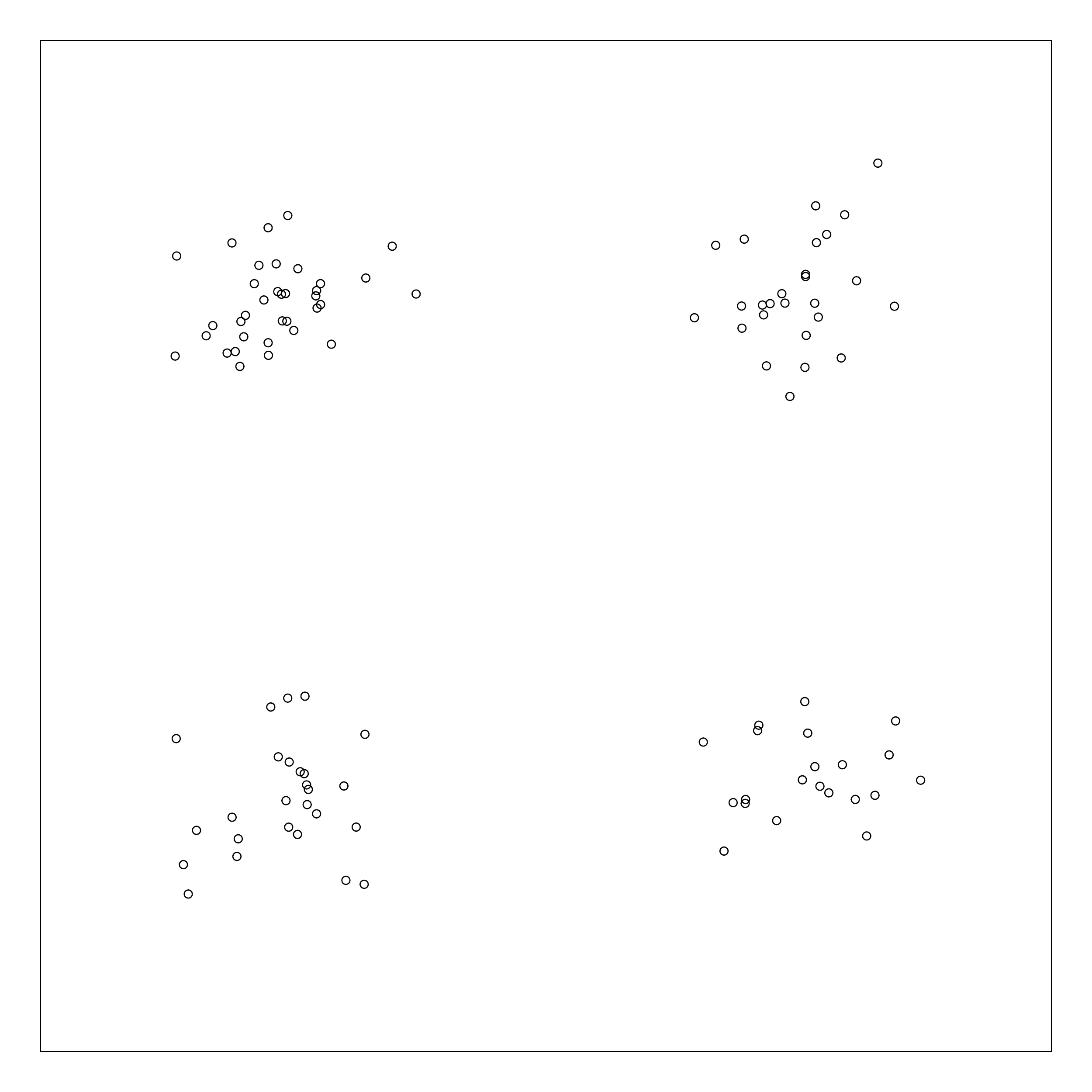}
  \includegraphics[scale=.14]{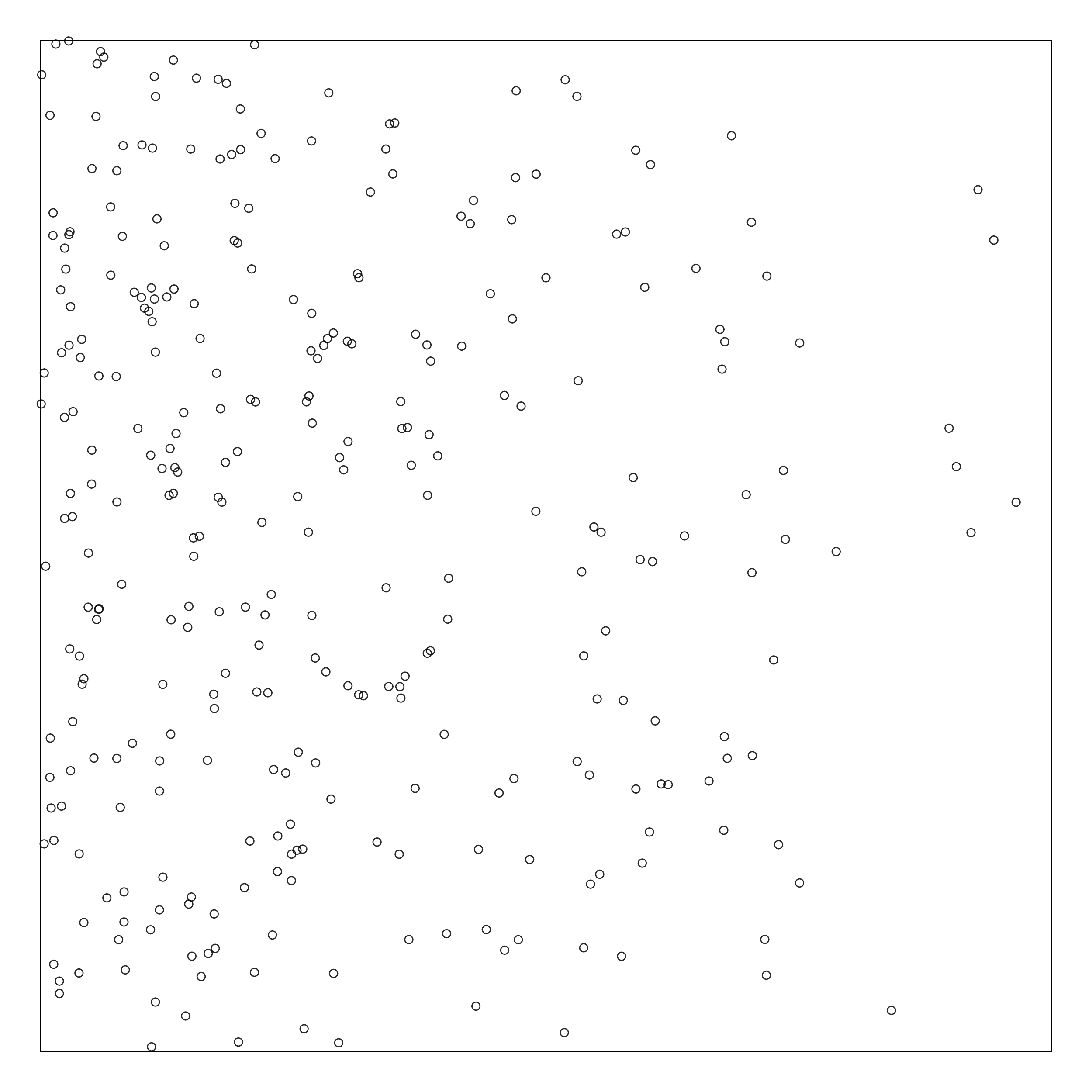}
  \caption{Realisations of a Poisson point process when its intensity function follows  Model 1 to Model 4 from left to right.}
  \label{fig_pois}
\end{figure}

We consider the non-parametric estimation of the intensity $\rho(x,y)$ of a spatial Poisson point process, for $(x,y)\in [0,1]^2$, given one realization of the process on $[0,1]^2$. Four models are considered: 
\begin{itemize}
\item Model 1: homogeneous with low intensity, $\rho(x,y)=100$.

\item Model 2: homogeneous with high intensity, $\rho(x,y)=1000$.

\item Model 3: four clusters.
Denoting by $\phi_{a,b}$ the bivariate Gaussian density centered at $(a,b)$ with standard deviation $0.05$, i.e. $\phi_{a,b}(x,y)= \exp(-((x-a)^2+(y-b)^2)/0.05^2)/(2\pi 0.05^2)$, 
$$\rho(x,y)=25(\phi_{0.25,0.25}(x,y)+\phi_{0.25,0.75}(x,y) + \phi_{0.75,0.25}(x,y)+\phi_{0.75,0.75}(x,y)).$$

\item Model 4: exponential decreasing on the $x$-axis,  $\rho(x,y)=1000  \exp(-3x)$.
\end{itemize}
Typical realizations of these 4 situations are shown in Figure~\ref{fig_pois}.

 \medskip
 
 The estimation of $\rho(x,y)$ is carried out using the standard kernel-based estimator (implemented in  \texttt{R}  by the function \texttt{density} of the package  \texttt{spatstat}), for which the choice of the bandwidth is crucial for the quality of estimation. 
 We consider three possibilities offered by  \texttt{spatstat} leading to the estimators $\hat\rho_1$, $\hat\rho_2$ and $\hat\rho_3$ respectively: The default one which is $1/8$ of the shortest length of the observation window, the choice \texttt{bw.diggle} suggested in \cite{diggle85}, and \texttt{bw.ppl}  based on likelihood cross-validation.
 
 In this functional estimation setting, the MSE matrix $\Sigma$ of the estimators is replaced by the MISE matrix, that we still denote by $\Sigma$,  with generic term $\Sigma_{ij} = \mathbb E \int(\hat \rho_i(x,y) - \rho(x,y))(\hat \rho_j(x,y) - \rho(x,y)){\rm d} x{\rm d} y$, $i,j=1,...,3$. 
 To average these 3 estimators, we estimate the optimal weights \eqref{lambda_one} using a bootstrap procedure to get $\hat\Sigma$, where we choose $\hat\rho_3$ as an initial estimator, as it is empirically the most robust. Specifically, given $\hat \rho_3$, $N$ independent samples of the Poisson point process with intensity $\hat\rho_3$ are simulated on the unit square ($N=100$ below), from which we deduce an estimation of the MISE matrix by discretizing the integral on the grid of estimation returned by \texttt{spatstat} (that is a 128x128 pixel array).
  
Given a sample $X$, the full procedure in \texttt{R} to get the final average estimator is as follows. It takes approximately 90 seconds on a regular laptop, for the considered models.
\begin{lstlisting}
#Computing the initial estimators
 est1 <-density(X) 
 est2 <-density(X,bw.diggle)
 est3 <-density(X,bw.ppl) 

#Bootstrapping the model to get a sample of the estimators
 N <-100
 ppboot <-rpoispp(est3,nsim=N)
 estboot1 <-lapply(ppboot,density)
 estboot2 <-lapply(ppboot,density,bw.diggle)
 estboot3 <-lapply(ppboot,density,bw.ppl)

#Deducing an estimation of the MISE matrix
 fun <-function(x,y,z){
   temp <-rbind(as.numeric(x-est3),as.numeric(y-est3),as.numeric(z-est3))
   return(temp %*% t(temp))}
 hatSigma <-matrix(0,3,3)
 for(i in 1:N){
   hatSigma <-hatSigma+fun(estboot1[[i]],estboot2[[i]],estboot3[[i]])}
  
#Constructing the average estimator
 invhatSigma <-solve(hatSigma)
 weights <-rowSums(invhatSigma)/sum(invhatSigma)
 AV <-weights[1]*est1+weights[2]*est2+weights[3]*est3
 AV[AV<0] <-0
\end{lstlisting}

Note that in the last step, the average estimator is projected on the space of positive functions to give a final consistent result in view of intensity estimation. The obtained estimator is therefore closest to the true intensity than the average estimator (by projection onto a convex set) and thus inherits its optimality properties. Repeating this procedure $10^3$  times for each model described above, we obtain an estimation of the MISE for each initial estimator and for the  (projected) average estimator, summarized in Table~\ref{sim:pois}. From this table, it appears that the initial estimators have variable performances, depending on the underlying intensity, but the most reliable one seems to be the choice of bandwidth based on likelihood cross-validation. The average estimator outperforms all initial estimators in all cases.

\begin{table}[h]
\centering
\begin{tabular}{|l|r|r|r|r|}
\hline
 & Raw & Diggle & PPL & AV \\ \hline
Model 1 & 729 (10.3) & 2903 (32) & 237 (11.8) &  \textbf{229} (9.7)  \\ \hline
Model 2 & 7240 (88) & 28673 (282) & 2356 (109) & \textbf{2247} (88)  \\ \hline
Model 3 & 54074 (38) & 14108 (199) & 12401 (111) &  \textbf{12081} (112)  \\ \hline
Model 4 & 115942 (347) &  137399 (482) &  116372 (408) &   \textbf{115762} (362)  \\ \hline
\end{tabular}
\caption{Estimation of the MISE for each initial estimator of the intensity  of a Poisson point process (given by Models 1-4), and for the average estimator ("AV"), based on $10^3$ replications. The initial estimators correspond to the kernel estimator for the bandwidth : "Raw" (default choice in the function \texttt{density} of \texttt{spatstat}), "Diggle" (option \texttt{bw.diggle}), "PPL" (option \texttt{bw.ppl}). An estimation of the standard deviation of the MISE estimation is given in parenthesis. }
\label{sim:pois}
\end{table}

\subsection{Determinantal point processes}

Determinantal point processes (DPPs) are models for regular point patterns. We refer to \cite{LMR15} for  their main statistical properties. A DPP is completely characterized by a kernel $C:\mathbb R^d\times \mathbb R^d\to \mathbb C$, which for existence needs to be a continuous covariance function and further satisfy a spectral condition. In the homogeneous case where $C(u,v)=C_0(v-u)$, this spectral condition reduces to $\mathcal F (C_0)\leq 1$ where $\mathcal F$ denotes the Fourier transform. In the non-homogeneous case, a sufficient condition is that there exists $C_0$ as before such that $C_0(u-v)-C(u,v)$ remains a covariance function. 

\begin{figure}
\centering
   \includegraphics[scale=.25]{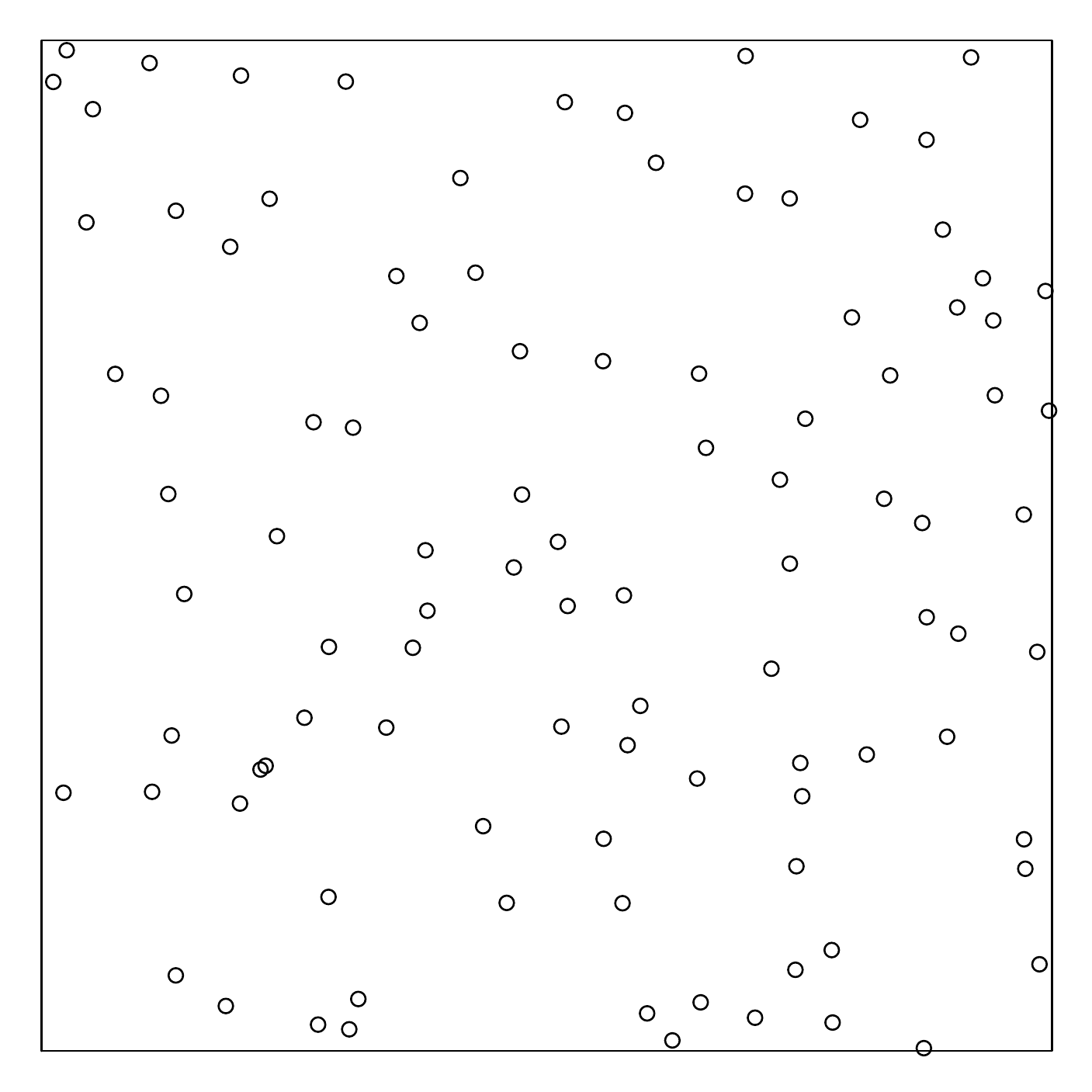}
  \includegraphics[scale=.25]{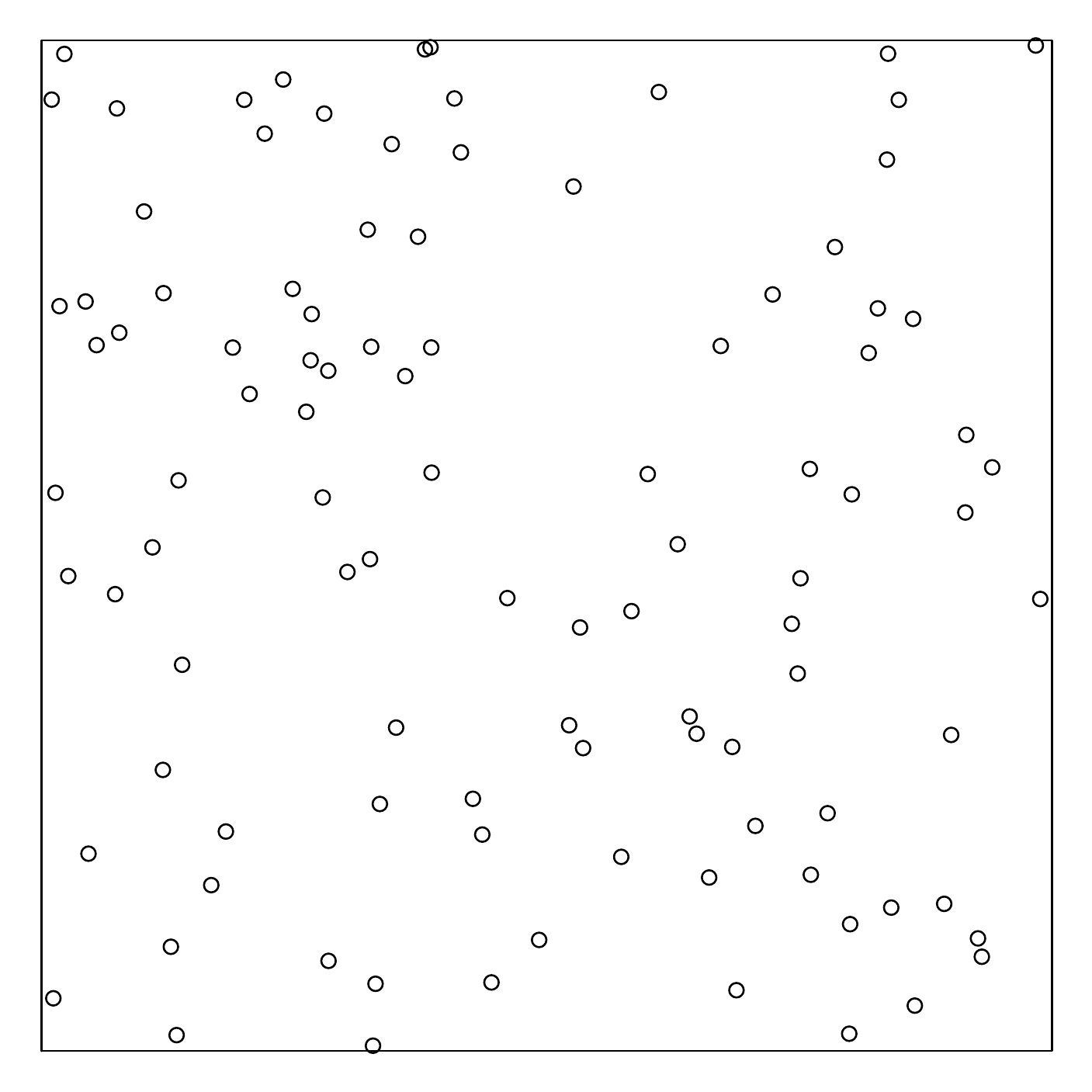} 
  \includegraphics[scale=.25]{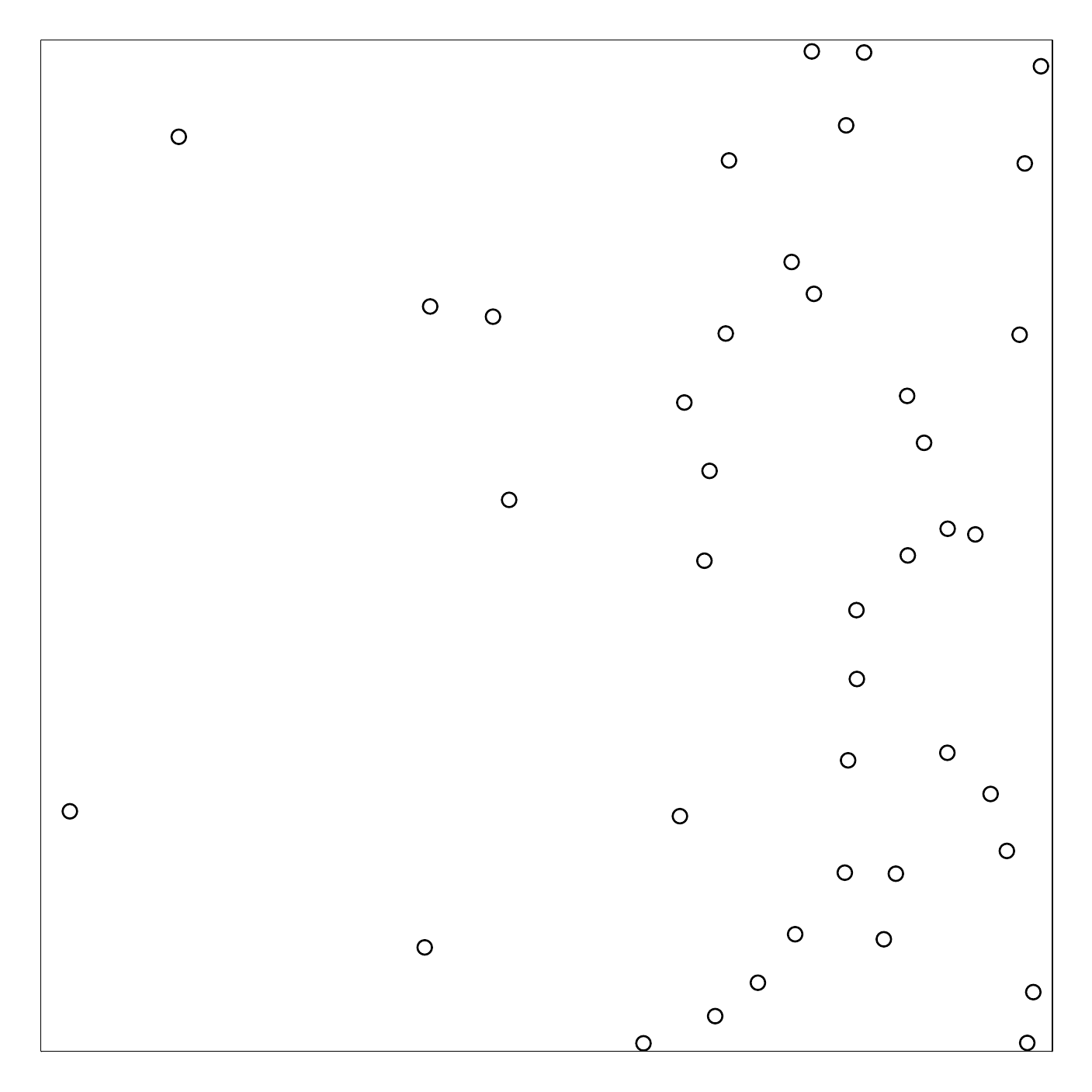}
  \includegraphics[scale=.25]{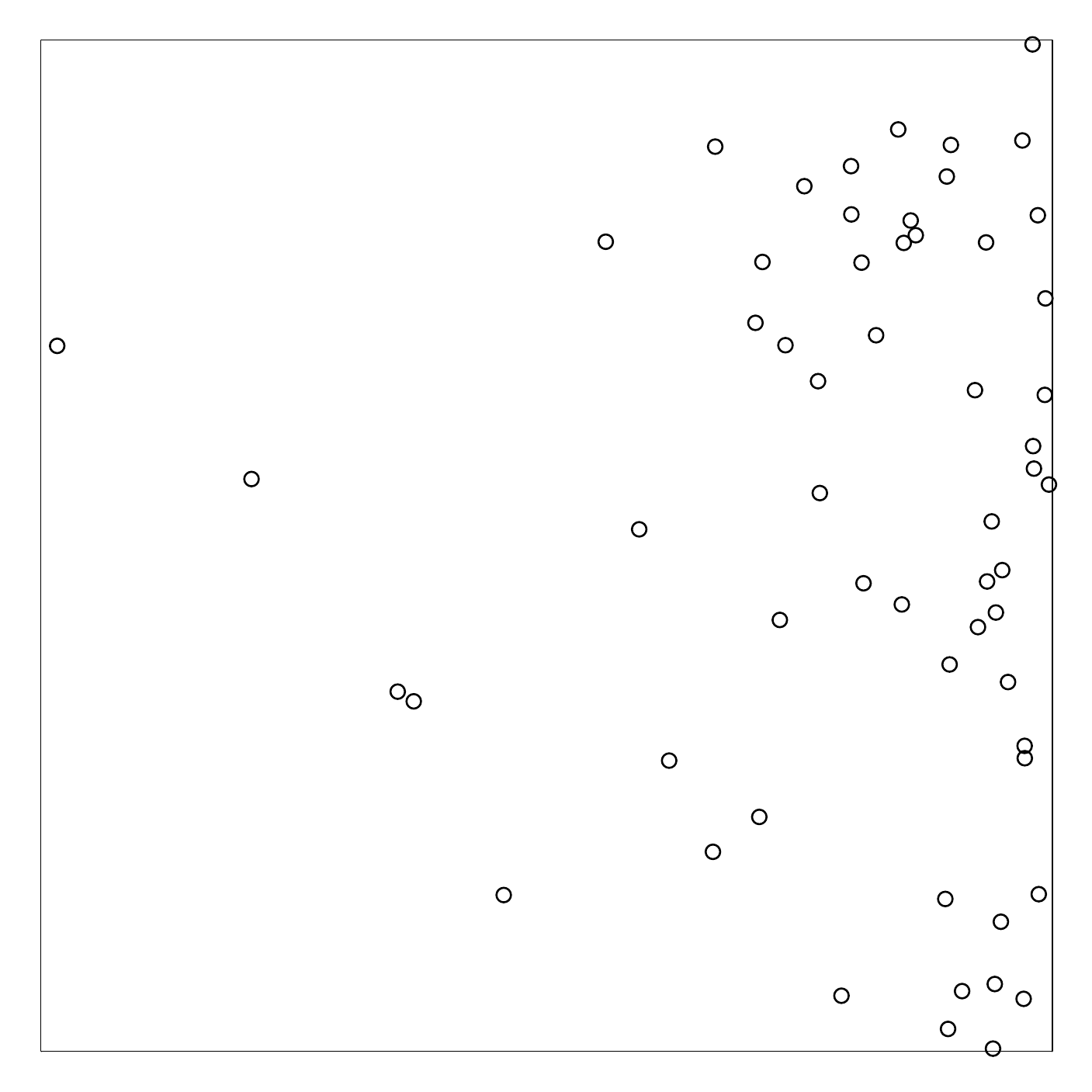}
  \caption{Realisations of DPP models on $[0,1]^2$ defined from left to right by DPP 1, DPP2, DPP3 and DPP4.}
\label{fig_DPP}
\end{figure}

In this section we consider the estimation of  parametric DPPs on the plane, defined through a parametric kernel  $C(u,v)=\sqrt{\rho(u)}\e^{-||u-v||^2/\alpha^2}\sqrt{\rho(v)}$   for $u,v\in\mathbb R^2$, where $\rho$ is assumed to be log-linear. For this model,  the intensity function is $\rho$ and the pair correlation function (pcf) is the isotropic function $g(r)=1-\e^{-2r^2/\alpha^2}$ for $r>0$, see \cite{LMR15}. Specifically, denoting $u=(x,y)$, we consider the four following situations.
 \begin{itemize}
\item DPP 1: $\rho(x,y)=100$ and $\alpha=\alpha^{(1)}_{\max}\approx 0.056$, which is an homogeneous DPP with the maximum possible value for the scale parameter $\alpha$ when $\rho=100$, deduced from  the spectral condition for existence discussed above.

\item  DPP 2: $\rho(x,y)=100$ and $\alpha=\alpha^{(1)}_{\max}/2$, which is the same model as above with less repulsion between the points. 

\item  DPP 3: $\rho(x,y)=4 \exp(4x)$ and $\alpha=\alpha^{(2)}_{\max}\approx 0.038$ which is an inhomogeneous DPP with  exponential increasing intensity along the $x$-axis, and the maximum possible value for the scale parameter $\alpha$.

\item  DPP 4: $\rho(x,y)=4 \exp(4x)$ and $\alpha=\alpha^{(2)}_{\max}/2$, the same model as DPP 3 with 
less repulsion between the points.
\end{itemize}
Typical realizations of these four processes on $[0,1]^2$ are shown in Figure~\ref{fig_DPP}.\\

In theory, the density of a DPP on any compact set is known, making possible likelihood estimation. However this density involves a new kernel, obtained from a spectral representation of $C$, which is rarely known in practice. In the homogeneous case, and when the domain of observation  is rectangular, some efficient approximations are introduced in \cite{LMR15} to make likelihood estimation feasible. In this situation, we recommend to use the maximum likelihood estimator. In the inhomogeneous case or when the observation window is not rectangular, likelihood inference seems difficult to implement, but alternative methods are available : minimum contrast estimation based on the Ripley's $K$ function, or on the pcf $g$, composite likelihood estimation, or Palm likelihood estimation. These methods are implemented in the function \texttt{dppm} of \texttt{spatstat}. While the first two methods have good theoretical backgrounds (see \cite{BL15-3}), the two others have not been justified yet from a theoretical perspective. From our experience, composite likelihood is not stable and we do not use it in the following. None of the three other methods is objectively better than the others and an averaging procedure makes sense.

Whatever the method,  the intensity function is estimated in \texttt{dppm} by maximizing the Poisson likelihood (that we assume to be log-linear), so that the three retained estimation methods differ only for the estimation of $\alpha$. We thus average these three estimations of $\alpha$ using \eqref{weight}, where  $\hat\Sigma$ is obtained by parametric bootstrap. In this procedure, we choose as initial parameter  of the model the output of the Palm likelihood estimation and we simulate $N=100$ samples. If $X$ denotes the observed point pattern, the associated code in $\texttt{R}$ is the following. This procedure takes approximately 3 minutes on a regular laptop.
\begin{lstlisting}
#Computing the three initial estimators
 fitg <- dppm(X~x, dppGauss, method="mincon", statistic="pcf", rmin=0.01, q=1/2)
 fitK <- dppm(X~x, dppGauss)
 fitpalm <- dppm(X~x, dppGauss, method="palm")

#Bootstrapping the initial estimators 
 N <-100
 ppboot <- simulate(fitg, nsim=N)
 fitbootK <- lapply(ppboot, function(y) dppm(y~x, dppGauss))
 fitbootg <- lapply(ppboot, function(y) dppm(y~x, dppGauss, method="mincon", statistic="pcf" , rmin=0.01, q=1/2))
 fitbootpalm <- lapply(ppboot, function(y) dppm(y~x, dppGauss,method="palm"))

#Deducing the MSE matrix
 alphafun <-function(x) x$fitted$fixedpar$alpha
 fitboot <-c(fitbootK,fitbootg,fitbootpalm)
 mat <-matrix(unlist(lapply(fitboot,alphafun))-alphafun(fitg),nrow=N)
 hatSigma <-t(mat)%*%mat/N

#Constructing the average estimator and its estimated MSE
 invhatSigma <-solve(hatSigma)
 weights <-rowSums(invhatSigma)/sum(invhatSigma)
 AV <-weights[1]*alphafun(fitK)+weights[2]*alphafun(fitg)+weights[3]*alphafun(fitpalm)
 MSE_AV <- 1/sum(invhatSigma)
\end{lstlisting}

In Table~\ref{tab:dpp} we have summarized the MSE of each estimator of the above procedure, based on $10^3$ replications. The performances of the initial estimators are variable, depending on the underlying model, but in all cases, the average estimator is better than the best initial estimator.

\begin{table}[h]
\centering
\begin{tabular}{|l|r|r|r|r|}
\hline
 & K & g & Palm & AG \\ \hline
DPP 1 &  3.29 (0.21) & 6.04 (0.37) & 2.56 (0.21)  & {\bf 2.20} (0.18)  \\ \hline
DPP 2  & 12.7 (0.54) & 8.32 (0.39) & 9.37 (0.45)  & {\bf 8.31} (0.39)  \\ \hline
DPP 3  & 19.1 (1.33) & 13.1 (1.02) & 7.52 (0.54)  & {\bf 6.91} (0.43)  \\ \hline
DPP 4 & 32.5 (0.56) & 27.3 (0.59) & 10.5 (0.45)  & {\bf 10.1} (0.45) \\ \hline
\end{tabular}
\caption{Estimated MSE for the estimation of the scale parameter $\alpha$ in DPPs models, based on $10^3$ replications. The estimators are the minimum contrast estimator based on $K$ ("K"), the one based on the pcf $g$ ("g"), the maximum Palm likelihood estimator ("Palm") and their  average ("AV") given by \eqref{weight}. An estimation of the standard deviation of the MSE estimation is given in parenthesis. Each entry has been multiplied by $10^5$ for ease of presentation. }
\label{tab:dpp}
\end{table}

\subsection{Thomas process}\label{sec:Thomas}

\begin{figure}[h]
\centering
  \includegraphics[scale=.17]{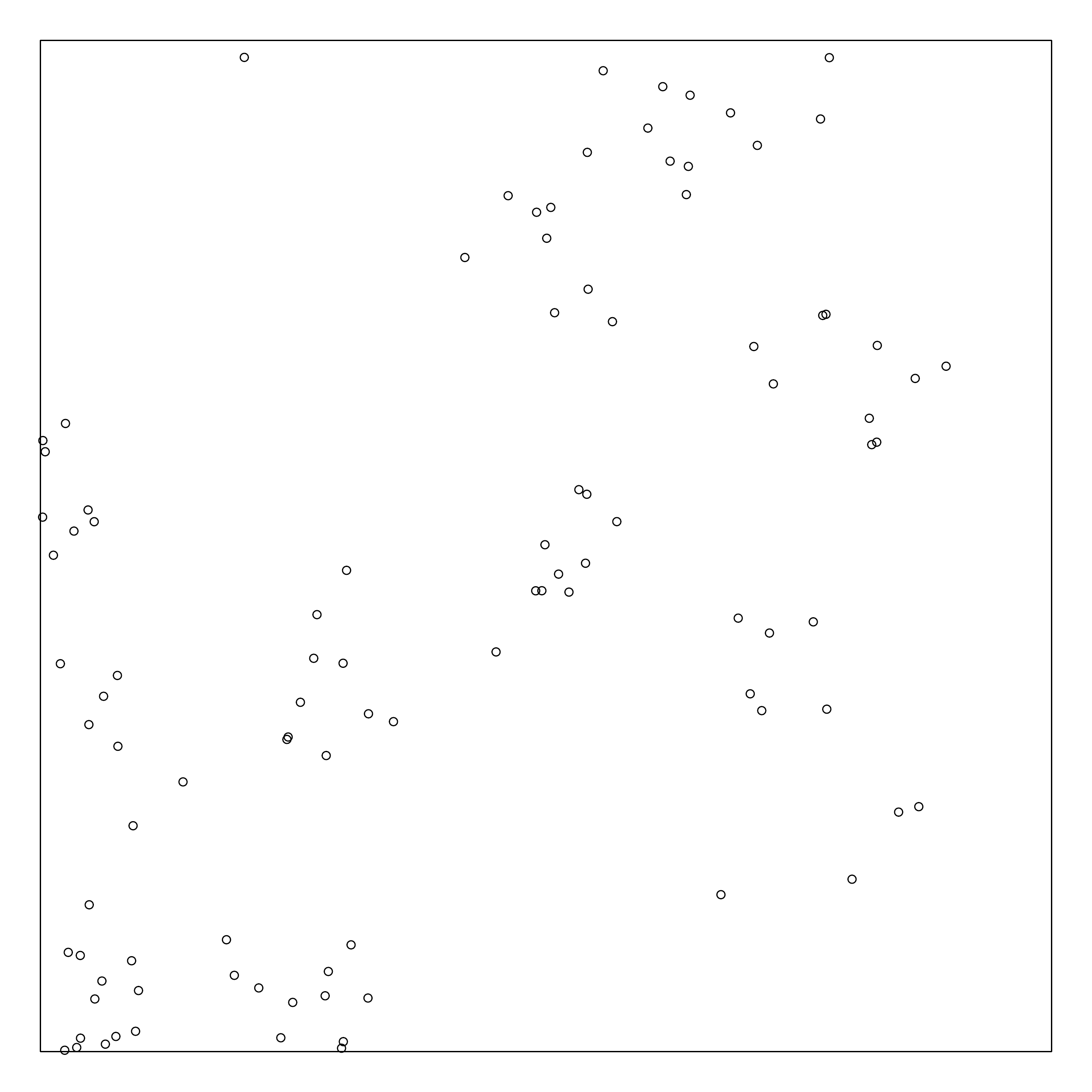}
  \includegraphics[scale=.17]{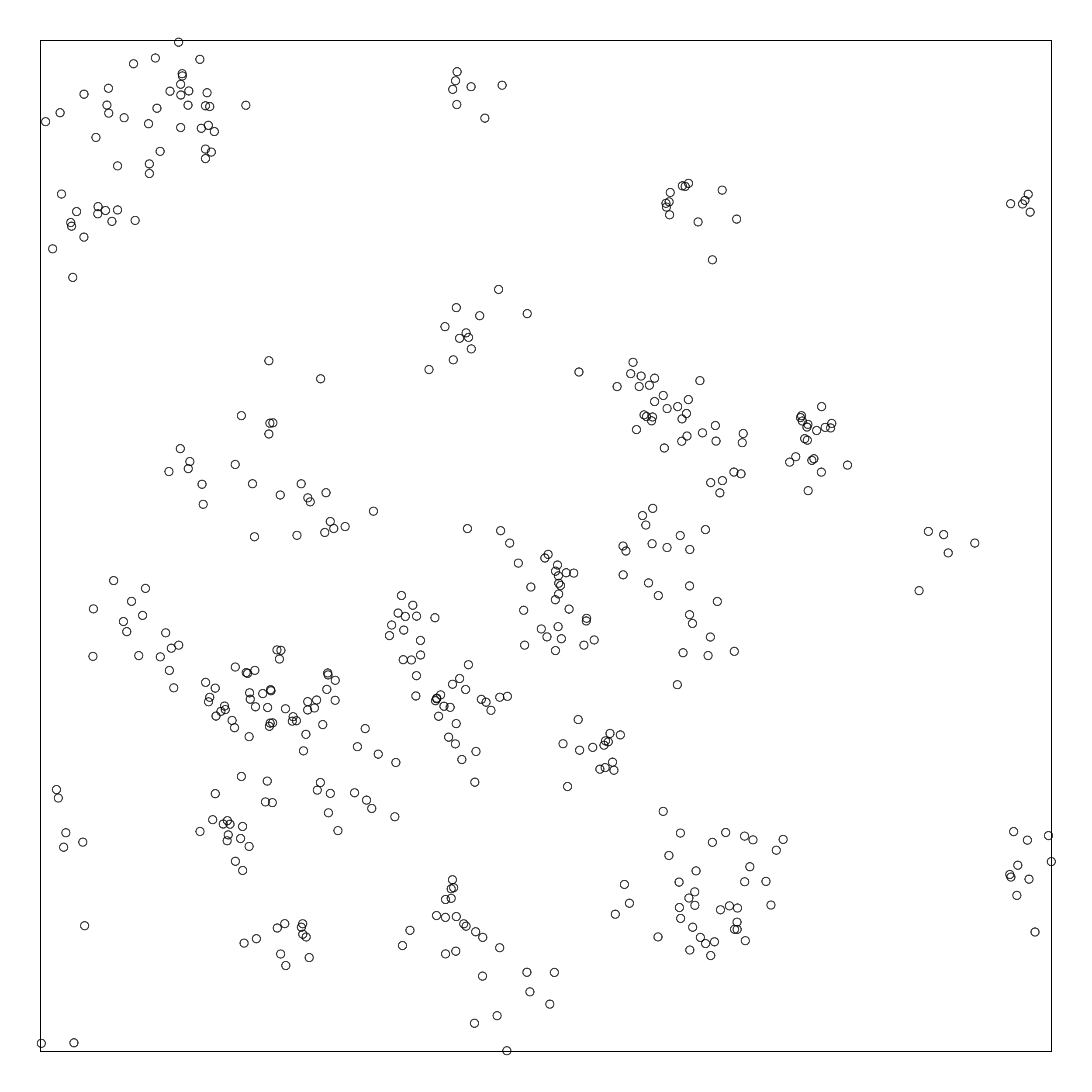}
  \includegraphics[scale=.17]{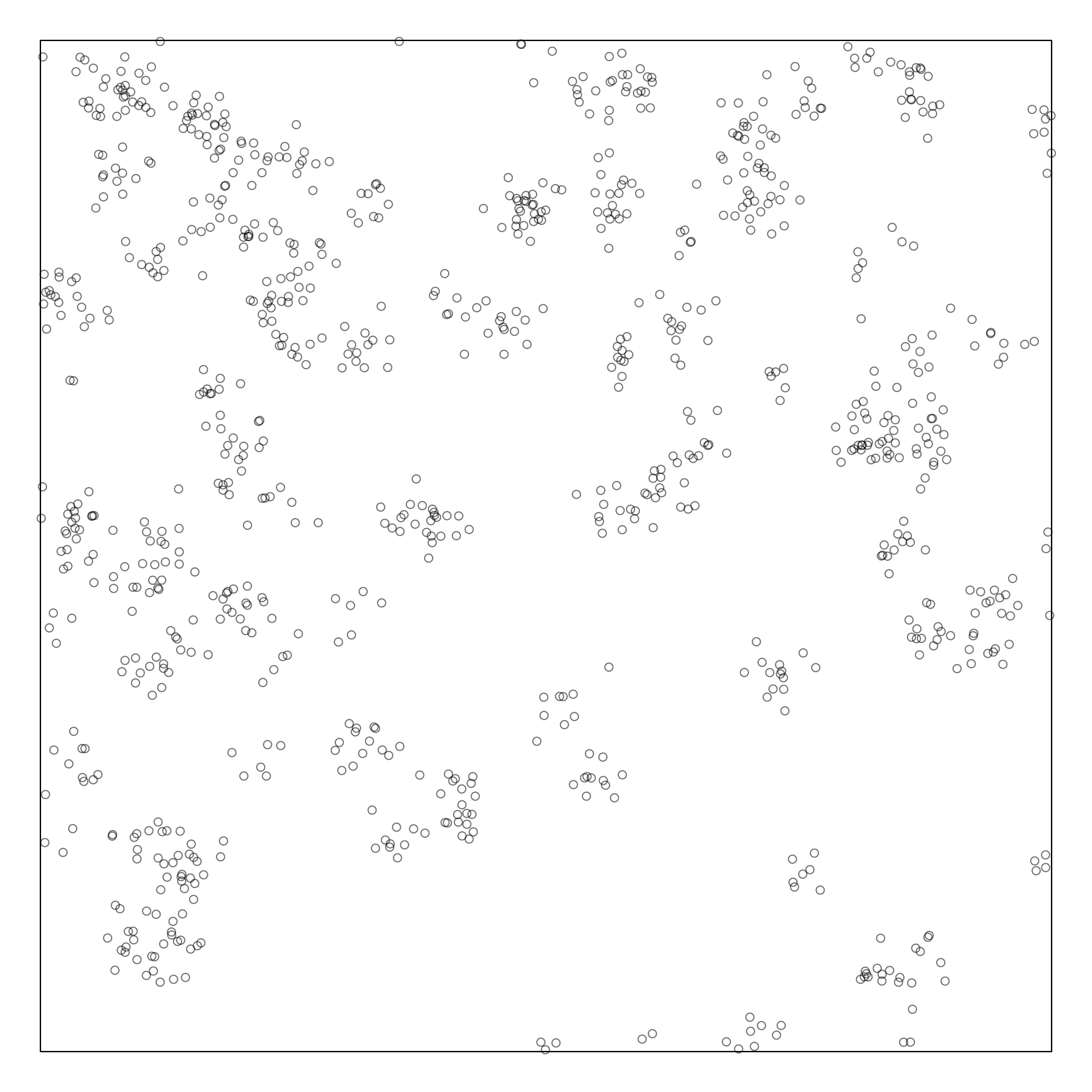}
  \caption{Realisations of a Thomas process with parameters $\kappa=10$, $\mu=10$ and $\sigma=0.05$ on, from left to right, $[0,1]^2$, $[0,2]^2$ and $[0,3]^2$.}
  \label{fig:Thomas}
\end{figure}

In this section, we consider the estimation of a Thomas process \cite{thomas}, which belongs to the larger class of Neyman-Scott  processes, see for instance \cite{moeller:waagepetersen:00}. This is a standard and classical example of model for clustered point patterns.
This model depends on three parameters: $\kappa$ represents the intensity of the "parents", generated as a homogeneous Poisson point process; $\mu$ is the mean number of points (or children) around each parent, drawn from a Poisson random variable; and $\sigma$ corresponds to the dispersion around each parent of his children. The children are sampled from a bivariate independent Gaussian distribution centered at the location of the parent with standard deviation $\sigma$. A realization of the Thomas process is given by the locations of the children, which are by construction organized by clusters.  A simulation is given in Figure~\ref{fig:Thomas} for $\kappa=10$, $\mu=10$, $\sigma=0.05$, and the observation window is $[0,1]^2$, $[0,2]^2$ and $[0,3]^2$, respectively.

Standard procedures to estimate the parameters of a Thomas process are minimum contrast estimation methods based  on $K$ or on  $g$, or maximum Palm likelihood estimation. These three methods are implemented in the function \texttt{kppm} of \texttt{spatstat}. Note that composite likelihood estimation is also proposed in this function, but from our experience the results are unstable and we do not use this method in the following. 

To average the above three estimators, we can either use \eqref{weight} for each parameter, or we can use the method described in Section~\ref{sec:foreign} to include all estimators for the estimation of each parameter, taking advantage of possible cross-correlations with the foreign estimators.  In each case, we decide to estimate the MSE matrix $\Sigma$ by parametric bootstrap where the initial estimator is the minimum contrast estimator based on the pcf $g$, and where we take $N=100$ samples.

Given a realization $X$ of a Thomas process, the procedure in \texttt{R} to get the average estimator accounting for foreign estimators is the following. It takes from 20 seconds for a point pattern as in the left hand side of Figure~\ref{fig:Thomas},  to two minutes for a point pattern as in the right hand side of the same figure.
\begin{lstlisting}
#Computing the initial estimators
 fitK <-kppm(X,~1,"Thomas",method="mincon",statistic="K")
 fitg <-kppm(X,~1,"Thomas",method="mincon",statistic="pcf")
 fitpalm <-kppm(X,~1,"Thomas",method="palm")
 
#Bootstrapping the model and the initial estimators
 N <-100
 ppboot <-simulate(fitg,nsim=N)
 fitbootK <-lapply(ppboot,kppm,~1,"Thomas",method="mincon",statistic="K")
 fitbootg <-lapply(ppboot,kppm,~1,"Thomas",method="mincon",statistic="pcf")
 fitbootpalm <-lapply(ppboot,kppm,~1,"Thomas",method="palm")
 
#Deducing an estimation of the MSE matrix
 kappadiff <-function(x) unlist(lapply(x,function(y) y$par[1])) - fitg$par[1]
 sigma2diff <-function(x) unlist(lapply(x,function(y) y$par[2])) - fitg$par[2] 
 mudiff <-function(x) unlist(lapply(x,function(y) y$mu)) - fitg$mu
 fitboot <-c(fitbootK,fitbootg,fitbootpalm)
 mat <-matrix(c(kappadiff(fitboot),sigma2diff(fitboot),mudiff(fitboot)),nrow=N)
 hatSigma <-t(mat)%*%mat/N
 
#Computing the full weights (taking into account foreign estimators)
 invhatSigma <-solve(hatSigma)
 matL <-kronecker(diag(1,3),rep(1,3))
 weights_full <-invhatSigma%*%matL%*%solve(t(matL)%*%invhatSigma%*%matL)
 
#Deducing the three average estimators and their estimated MSE
 param <-function(x) unlist(parameters(x))[-1]
 estvec <-as.vector(t(sapply(list(fitg,fitK,fitpalm),param)))
 AV_plus <-t(weights_full)%*%estvec
 MSE_AV_plus <-t(weights_full)%*%hatSigma%*%weights_full
 \end{lstlisting}
To get the average estimators that do not use foreign estimators, the above code differs only in the  last two steps:
\begin{lstlisting}
#Computing the componentwise weights (without foreign estimators)
 support=kronecker(diag(1,3),matrix(1,3,3))
 hatSigma_sparse <-hatSigma*support
 invhatSigma_sparse <-solve(hatSigma_sparse)
 weights_sparse <-invhatSigma_sparse%*%matL%*%solve(t(matL)%*%invhatSigma_sparse%*%matL)


#Constructing the average estimators and its estimated MSE
 param <-function(x) unlist(parameters(x))[-1]
 estvec <-as.vector(t(sapply(list(fitg,fitK,fitpalm),param)))
 AV <-t(weights_sparse)%*%estvec
 MSE_AV <-t(weights_sparse)%*%hatSigma%*%weights_sparse
 \end{lstlisting}

In Table~\ref{tab:Thomas}, we give the estimated MSE of each initial estimator and  of the average estimator whether it uses foreign estimators (AV+) or not (AV). This table is based on $10^3$ replications of the above procedure when the observation window is 
either $[0,1]^2$, or $[0,2]^2$ or $[0,3]^2$. From this study, it turns out that the contrast estimation method based on $g$ is the best among the three initial estimators, but it is globally outperformed by the average estimator, with or without the inclusion of foreign estimators. Further remarks are in order. While AV+ should in theory (i.e. if $\Sigma$ were perfectly known) be better than AV, this is not necessarily the case in practice. There are two possible reasons. The first situation occurs if there is not enough data to hope for a good estimation of the full MSE matrix $\Sigma$. In our example, this matrix  contains 45 unknown quantities and  its estimation may clearly be inaccurate for small data sets, as when $L=1$ in Table~\ref{tab:Thomas}. A second reason is when AV+ is in fact more or less equal to AV in theory, meaning that the weights associated to the foreign estimators should be zero (this is for instance the case if there are no correlations with the foreign estimators). In this situation, the inclusion of foreign estimators can be viewed as a noise in the averaging procedure that can only deteriorate the estimation. This is what happens for the estimation of $\kappa$, when $L=2$ and $L=3$ in Table~\ref{tab:Thomas}. However, when the data are rich enough, the weights are sufficiently well estimated so that AV+ can be expected to be at least as good as AV.

\begin{table}[h]
\centering
\begin{tabular}{|ll|r|r|r|r|r|}
\hline
& & $K$ & $g$ & Palm & AV & AV+ \\ \hline
$L=1$ & $\kappa$ & 51.66 (4.95) & 45.39 (3.86) &  52.90 (4.48) & 40.49 (4.03) & \textbf{34.41} (2.56) \\ 
 & $\sigma^2$ & 11.87 (0.76) & \textbf{11.57} (0.81) & 19.83 (2.23) &  12.79 (3.34) & 16.24 (1.22)  \\ 
 & $\mu$ &  19.62 (1.28) & 20.00 (1.43) & 20.68 (1.39) &  \textbf{19.16} (1.27) & 19.55 (1.25)\\ \hline
 $L=2$ & $\kappa$ &  15.13 (0.98) &  9.12 (0.53) &  12.27 (0.71) &  \textbf{7.60} (0.46) & 9.18 (0.57) \\ 
  & $\sigma^2$  & 6.49 (0.41) & 2.84 (0.12) & 7.38 (2.10) &  2.75 (0.11) & \textbf{2.29} (0.10) \\ 
   & $\mu$ & 8.00 (0.40) & 5.40 (0.30) &  6.29 (0.37) &  5.62  (0.35) & \textbf{4.87} (0.24) \\ \hline
 $L=3$ & $\kappa$ &  6.57 (0.34) & 3.31 (0.15) & 5.50 (0.26) &  \textbf{3.06} (0.14) & 3.39 (0.16) \\ 
  & $\sigma^2$  & 4.90 (0.41) & 1.40 (0.10) & 2.54 (0.24) &  1.17 (0.07) & \textbf{1.04} (0.06) \\ 
   & $\mu$ &  4.87 (0.28) & 2.52 (0.14) & 3.07 (0.14) &  2.47 (0.14) & \textbf{2.18} (0.10) \\ \hline
\end{tabular}
\caption{Estimated MSE of the estimators of the parameters of a Thomas process, observed on $[0,L]^2$, based on $10^3$ replications. The average estimator  includes foreign estimators (AV+) or not (AV). An estimation of the standard deviation of the MSE estimation is given in parenthesis. The entries for $\sigma^2$ has been multiplied by $10^7$ for ease of presentation. }
\label{tab:Thomas}
\end{table}

In conclusion, for the Thomas process, we recommend to use the standard averaging procedure AV for each parameter (without foreign estimators), given by \eqref{weight}, which safely provides a better estimate in all situations. If the data are rich enough (like in the case of the observation window $[0,3]^2$ in our example), the inclusion of foreign estimators seems reasonable and may improve the final result.

\subsection{Boolean model}\label{sec:bool}

The following simulation study  is already presented in \cite{lavancier2016general}, but we judged interesting to include it in the present contribution, since it concerns the main model of random sets used in spatial statistics and stochastic geometry, see \cite{chiu2013}.  It is moreover a good example of situation where the inclusion of foreign estimators in the averaging procedure is highly relevant.

\begin{figure}[!htbp]%
  \centering 
  
   \includegraphics[scale=.21]{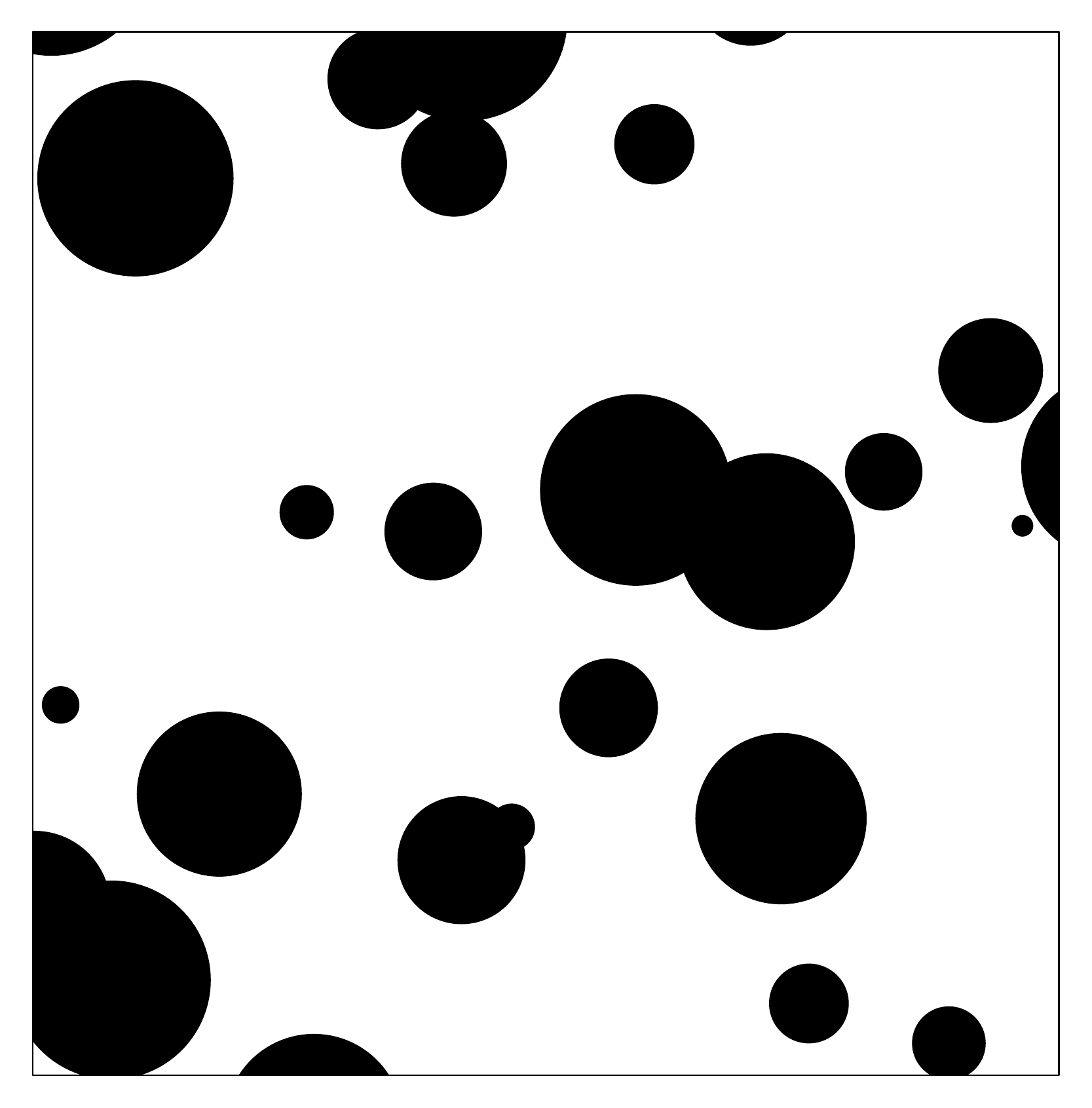}
    \includegraphics[scale=.21]{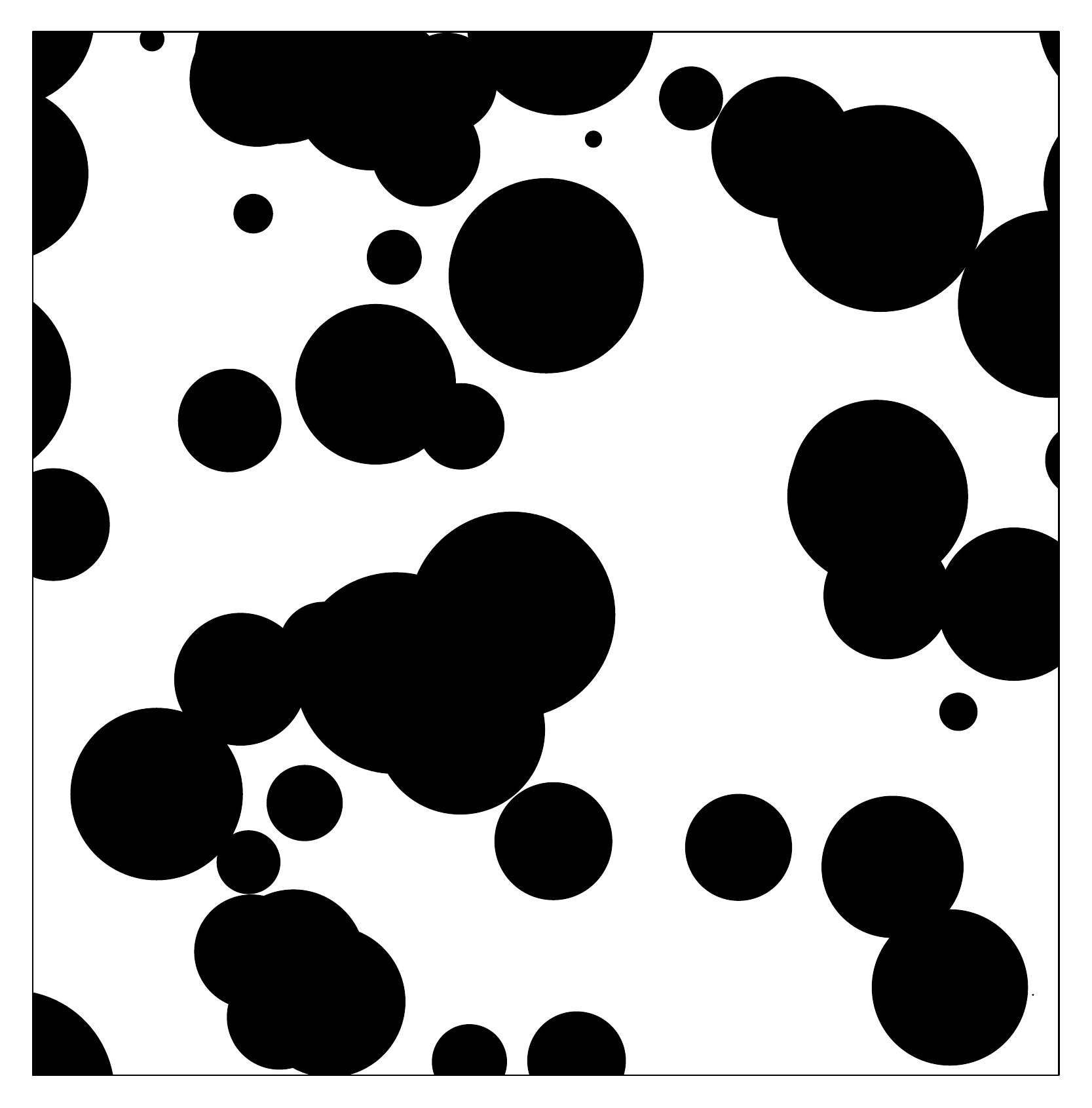}
       \includegraphics[scale=.21]{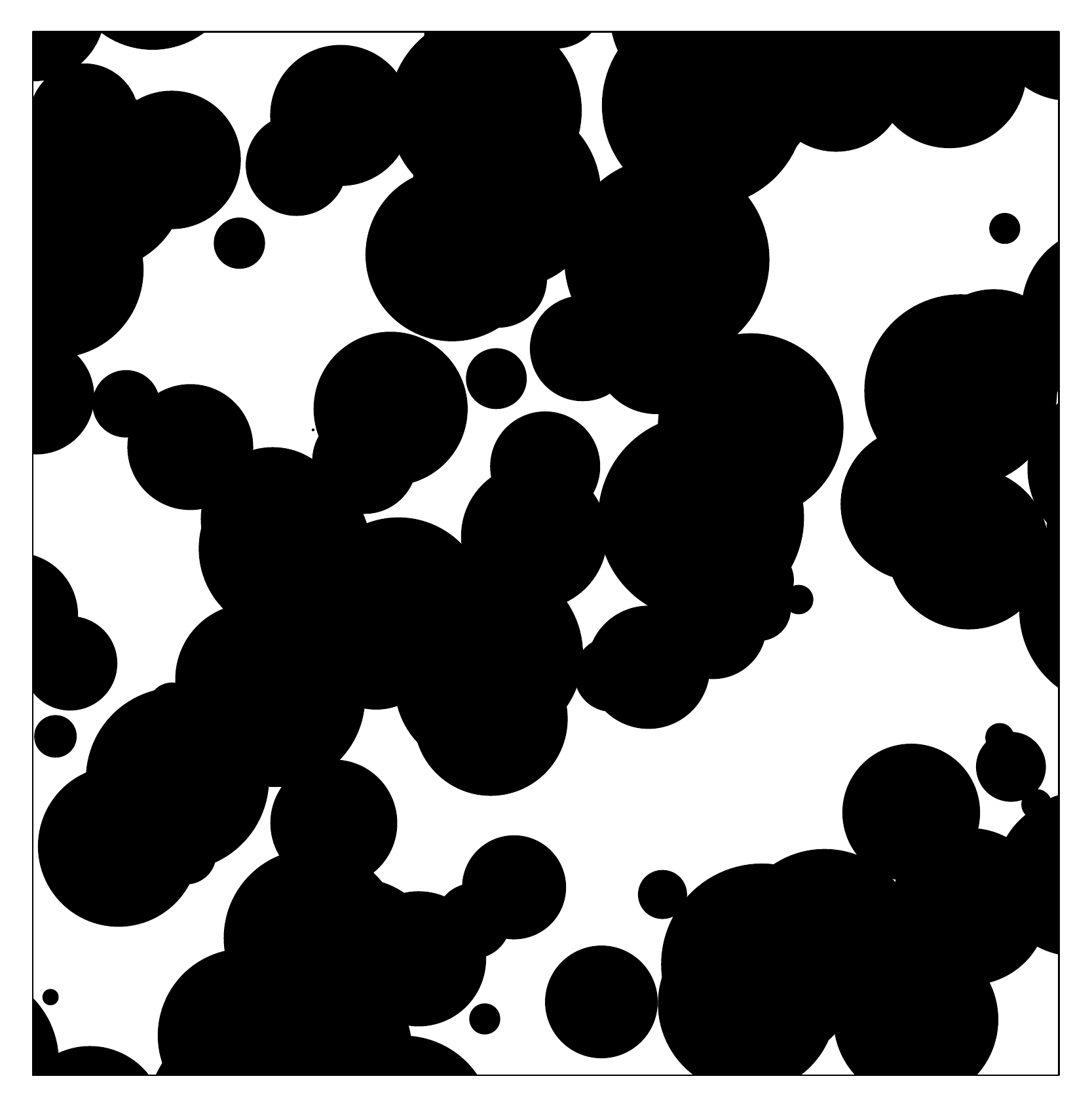}
    \includegraphics[scale=.21]{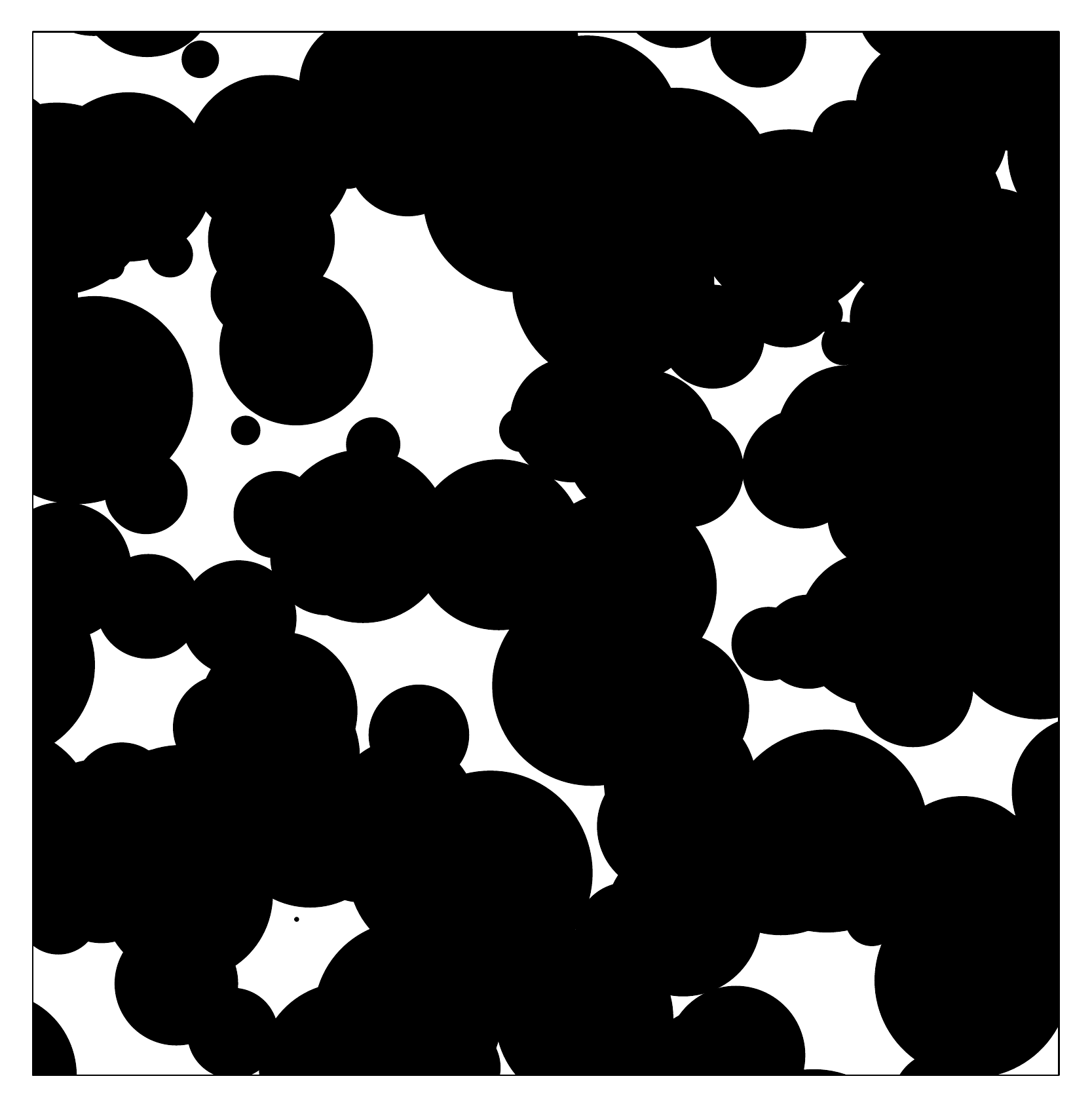}

  \caption{Samples from a Boolean model on $[0,1]^2$ with intensity, from left to right,  $\rho=25, 50, 100, 150$ and law of radii $B(1,\alpha)$ where $\alpha=1$.}
  \label{fig:bool}
\end{figure}

We consider a planar Boolean model where the germs come from a homogeneous Poisson point process with intensity $\rho$ and the grains are independent random discs, the radii of which are distributed according to a beta distribution over $[0,0.1]$ with parameter $(1,\alpha)$, $\alpha>0$, denoted  by $B(1,\alpha)$.  Figure~\ref{fig:bool} contains four realizations of this model on $[0,1]^2$ where $\rho=25, 50, 100, 150$ respectively and $\alpha=1$.\\

%
%

To estimate $\rho$ and $\alpha$,  we use a moment method based on the perimeter and the area of the random set generated by the Boolean model,  see \cite{Molchanov97}, and we denote the result by $\hat\rho_1$ and $\hat\alpha$.  As an alternative procedure to estimate $\rho$, it is also possible to apply  a method  introduced in  \cite{molchanov1995}, based on the number of tangent points of the random set in certain directions, which provides a new estimator $\hat\rho_2$. These estimation methods are detailed for the present model in \cite{lavancier2016general}. 
We have finally two estimators $\hat\rho_1$ and $\hat\rho_2$ for $\rho$, and one estimator $\hat\alpha$ for $\alpha$.

The averaging method allows us to combine the two estimators of $\rho$ to improve the estimation of this parameter. Note that the inclusion of the foreign estimator $\hat\alpha$ for this purpose is not possible, since the sum of weights for foreign estimators must be zero, therefore in presence of only one foreign estimator, its weight is zero. On the other hand, it is also possible to improve the estimation of $\alpha$, even if there is only one available estimator for it, by including the two foreign estimators $\hat\rho_1$ and $\hat\rho_2$. In this case the average estimator of  $\alpha$ is of the form $\hat\alpha_{AV}=\hat\alpha+ \mu (\hat\rho_1 - \hat\rho_2)$ where $\mu$ is the weight to be estimated. The full averaging procedure thus corresponds to the framework of Section~\ref{sec:foreign} where $J=2$ and $K=1$, in particular the optimal weights are given by \eqref{oracle_max}.
To estimate $\Sigma$, we use a parametric bootstrap of the model using $N=100$ samples with initial parameters $0.5 (\hat\rho_1+\hat\rho_2)$ and $\hat\alpha$.  Concerning  the implementation of this procedure, the main task is to code the initial estimators, which is not a procedure available by default on \texttt{R}. We do not enter into these details here. Then it is straightforward to deduce $\hat\Sigma$ as in the previous sections and to get the average estimator as detailed in Section~\ref{sec:foreign}.

Table~\ref{tab:bool} reports the MSE of each estimator, estimated from $10^4$ replications of a Boolean model on $[0,1]^2$ with parameters $\rho=25, 50, 100, 150$ and $\alpha=1$. 
From these results, the average estimators clearly outperform the initial estimators. The improvement is in particular remarkable for high values of $\rho$, i.e. in presence of a dense random set. 
 The parameter $\alpha$ ruling the law of the radii is known to be difficult to estimate, in particular when $\rho$ is high. The fact that $\rho$ can be easier to estimate together with the cross-correlations between $\hat\alpha$, $\hat\rho_1$ and $\hat\rho_2$ make relevant the inclusion of foreign estimators in $\hat\alpha_{AV}$. As demonstrated in  Table~\ref{tab:bool}, the difference of efficiency between $\hat\alpha$ and $\hat\alpha_{AV}$ can be impressive.

\begin{table}[htbp]
\begin{center}
\begin{tabular}{|l|ccc|cc|}
\hline
  & $\hat \rho_1$ & $\hat \rho_2$ & $\hat \rho_{AV}$ & $\hat\alpha$  &  $\hat\alpha_{AV}$ \\ \hline
\multirow{2}{*}{$\rho=25$}  &  34.15 & 14.63 & {\bf 14.60} & 8.09 & {\bf 6.70} \\ 
 & (0.55) & (0.22) & (0.22) & (0.15) & (0.13)  \\ \hline
\multirow{2}{*}{$\rho=50$}  &  131.63 & 47.41 & {\bf 45.65} & 4.69 & {\bf 3.24}  \\ 
 & (2.26) & (0.72) & (0.67) & (0.067) & (0.048) \\ \hline
\multirow{2}{*}{$\rho=100$} &  949 & 272 & {\bf 223 } & 5.70 &  {\bf 2.29}  \\ 
 & (21.8) & (4.9) & (3.6) & (0.086) & (0.034) \\ \hline
\multirow{2}{*}{$\rho=150$} &  7606 & 1656 & {\bf 1005} & 14.7 &  {\bf 4.1} \\ 
 & (341) & (46.5) & (24.4) & (0.34) & (0.11)  \\ \hline
\end{tabular}
\end{center}
\caption{Estimated MSE of the initial estimators $\hat\rho_1$, $\hat\rho_2$, $\hat\alpha$ and of the average estimators $\hat\rho_{AV}$ and $\hat\alpha_{AV}$ based on  $10^4$ replications of a Boolean model with intensity  $\rho=25,\, 50,\, 100,\, 200$ and law of radii $B (1,\alpha)$ with $\alpha=1$. An estimation of the standard deviation of the MSE estimation is given in parenthesis.  The last two columns have been multiplied by 100 for ease of presentation.}
\label{tab:bool}
\end{table}

\section{Conclusion}\label{sec:conclusion}
The objective of averaging is to produce a single final efficient estimator in a statistical inference problem for which several methods are available. The solution is constructed from an estimation of the mean-square error matrix  of the initial estimators. In most cases, the average estimator improves on the best estimator in the collection, making averaging suitable even in situations where one of the initial estimators is known to be better than the rest. When implemented carefully, it is rare to see the averaging procedure perform truly worse than the best initial estimator. Not the least, as a free by-product of the procedure, an estimation of the mean square error of the average estimator is available, see~\eqref{MSEAV}, making straightforward the construction of confidence intervals.\\

The computational cost relies almost entirely on the ability to estimate the MSE matrix, which varies according to the model. In a parametric model where $\Sigma$ can be expressed as a function of the parameters, this estimation reduces to plug-in, which is  generally easy and fast to compute. In most parametric models though, as those presented in Section~\ref{sec:appli}, the latter function is not explicit and the estimation of $\Sigma$  requires the use of re-sampling, i.e. parametric bootstrap. In all cases, the computational cost is comparable to the cost inherent to the estimation of the variance of an estimator, or the construction of confidence intervals.  For the same price, the averaging procedure allows to get a more accurate estimate together with an estimation of its MSE.\\

For the standard models of spatial statistics considered in Section~\ref{sec:appli}, the averaging procedure works well, as demonstrated in our simulation study. In general, an ideal situation to apply the averaging procedure should fulfill the two important requirements listed below. 
\begin{enumerate}
	\item \textit{The oracle improves over each initial estimator.} It goes without saying that the averaging process is only interesting if the objective (the oracle) can improve on the current methods. If this condition is not verified, a method to select the best estimator, less sensible to the estimation of $\Sigma$, is generally more appropriate. While the relative performance of the oracle may depend on some unknown factors, it can be expected to rely for the most part on the statistical model itself so that it can be investigated independently from the data.

	\item \textit{The mean-square errors of the estimators is estimable.} The accuracy of $\hat \Sigma$ is arguably the main factor for the efficiency of the averaging procedure. If the data do not enable to build a proper estimate of $\Sigma$,  seeking for a suitable combination is usually hopeless. 
This requirement is not as strong as it may seem.  If the MSE of the oracle is much lower than the MSE of the initial estimators (see the first point above), given that \eqref{mse} is very smooth in $\lambda$, it is expected that the MSE associated to the weight $\hat\lambda$ still remains significantly lower than the MSE of the initial estimators, even if $\hat\lambda$ is not so close to $\lambda^*$. 
Nonetheless, a better estimated $\Sigma$ unequivocally leads to a better average estimator. In particular, in presence of many initial estimators (thus inducing a large matrix $\Sigma$), we recommend to perform a pre-selection or to use a more sophisticated averaging procedure as discussed in Section~\ref{sec:complement}, unless there is a large amount of data. The same recommendation applies for the introduction, or not,  of foreign estimators, which increases the size of $\Sigma$. 
\end{enumerate}

These two conditions are essentially inherent to the statistical model at hand and can be investigated independently from the data in most situations. When applied to a new specific statistical model, we therefore recommend to perform a preliminary  analysis, which can be theoretical and/or involve a computational study,  to establish if these requirements are verified. This analysis may also serve to calibrate the averaging procedure, whether it concerns the choice of the constraint set of weights, the best way to estimate $\Sigma$ (in particular the choice of the plug-in estimator in a parametric bootstrap procedure), or the use of foreign estimators.

\bibliographystyle{acm}
\bibliography{biblio_spatial}

\begin{thebibliography}{10}

\bibitem{baddeley:rubak:turner:15}
{\sc Baddeley, A., Rubak, E., and Turner, R.}
\newblock {\em Spatial Point Patterns: Methodology and Applications with R.}
\newblock Chapman and Hall/CRC Press, London, 2015.

\bibitem{A-BadTur05}
{\sc Baddeley, A., and Turner, R.}
\newblock {M}odelling spatial point patterns in \texttt{R}.
\newblock {\em Journal of Statistical Software 12}, 6 (2005), 1--42.

\bibitem{bates1969combination}
{\sc Bates, J.~M., and Granger, C.~W.}
\newblock The combination of forecasts.
\newblock {\em Operations Research 20}, 4 (1969), 451--468.

\bibitem{BL15-3}
{\sc Biscio, C. A.~N., and Lavancier, F.}
\newblock Contrast estimation for parametric stationary determinantal point
  processes.
\newblock {\em to appear in Scandinavian Journal of Statistics
  (arXiv:1510.04222)\/} (2016).

\bibitem{MR2351101}
{\sc Bunea, F., Tsybakov, A.~B., and Wegkamp, M.~H.}
\newblock Aggregation for {G}aussian regression.
\newblock {\em Ann. Statist. 35}, 4 (2007), 1674--1697.

\bibitem{chiu2013}
{\sc Chiu, S.~N., Stoyan, D., Kendall, W.~S., and Mecke, J.}
\newblock {\em Stochastic geometry and its applications}, 3~ed.
\newblock John Wiley \& Sons, 2013.

\bibitem{diggle85}
{\sc Diggle, P.}
\newblock A kernel method for smoothing point process data.
\newblock {\em Applied Statistics 34}, 2 (1985), 138--147.

\bibitem{elliott2011averaging}
{\sc Elliott, G.}
\newblock Averaging and the optimal combination of forecasts.
\newblock Tech. rep., UCSD Working Paper, 2011.

\bibitem{gaiffas}
{\sc Ga{\"\i}ffas, S., and Lecu{\'e}, G.}
\newblock Hyper-sparse optimal aggregation.
\newblock {\em Journal of Machine Learning Research 12\/} (2011), 1813--1833.

\bibitem{hansen2007least}
{\sc Hansen, B.~E.}
\newblock Least squares model averaging.
\newblock {\em Econometrica 75}, 4 (2007), 1175--1189.

\bibitem{hjort2003frequentist}
{\sc Hjort, N.~L., and Claeskens, G.}
\newblock Frequentist model average estimators.
\newblock {\em Journal of the American Statistical Association 98}, 464 (2003),
  879--899.

\bibitem{LMR15}
{\sc Lavancier, F., M{\o}ller, J., and Rubak, E.}
\newblock Determinantal point process models and statistical inference.
\newblock {\em Journal of the Royal Statistical Society, series B 77}, 4
  (2015), 853--877.

\bibitem{lavancier2016general}
{\sc Lavancier, F., and Rochet, P.}
\newblock A general procedure to combine estimators.
\newblock {\em Computational Statistics \& Data Analysis 94\/} (2016),
  175--192.

\bibitem{molchanov1995}
{\sc Molchanov, I.~S.}
\newblock Statistics of the boolean model: from the estimation of means to the
  estimation of distributions.
\newblock {\em Advances in applied probability\/} (1995), 63--86.

\bibitem{Molchanov97}
{\sc Molchanov, I.~S.}
\newblock {\em Statistics of the Boolean Model for Practitioners and
  Mathematicians}.
\newblock Wiley, Chichester, 1997.

\bibitem{moeller:waagepetersen:00}
{\sc M{\o}ller, J., and Waagepetersen, R.~P.}
\newblock {\em Statistical Inference and Simulation for Spatial Point
  Processes}.
\newblock Chapman and Hall/CRC, Boca Raton, 2004.

\bibitem{R}
{\sc {R Core Team}}.
\newblock {\em R: A Language and Environment for Statistical Computing}.
\newblock R Foundation for Statistical Computing, Vienna, Austria, 2016.

\bibitem{quadprog}
{\sc {S original by Berwin A. Turlach R port by Andreas Weingessel}}.
\newblock {\em quadprog: Functions to solve Quadratic Programming Problems.},
  2013.
\newblock R package version 1.5-5.

\bibitem{thomas}
{\sc Thomas, M.}
\newblock A generalization of {P}oisson's binomial limit for use in ecology.
\newblock {\em Biometrika 36}, 1/2 (1949), 18--25.

\bibitem{timmermann2006forecast}
{\sc Timmermann, A.}
\newblock Forecast combinations.
\newblock In {\em Handbook of Economic Forecasting}, G.~Elliott, C.~Granger,
  and A.~Timmermann, Eds. North Holland, Amsterdam, 2006, pp.~135--196.

\bibitem{yang2004aggregating}
{\sc Yang, Y.}
\newblock Aggregating regression procedures to improve performance.
\newblock {\em Bernoulli 10}, 1 (2004), 25--47.

\end{thebibliography}

\end{document}